\documentclass[a4paper, 11pt]{article}
\usepackage{a4}
\usepackage{tikz}
\usepackage{amssymb}
\usepackage{amsmath}
\usepackage{epsfig}
\usepackage[english]{babel}
\usepackage{graphicx}
\usepackage{float}
\usepackage{algorithmic}
\usepackage{algorithm}

\usepackage[utf8]{inputenc} 
\usepackage{dsfont} 
\usepackage{color} 
\usepackage{authblk}
\usepackage{bm}
\usepackage{subfig}
\usepackage{amsmath}
\usepackage{verbatim}

\newcommand{\Var}{\text{Var}}
\newcommand{\Cov}{\text{Cov}}

\newcommand{\btabc}{\begin{tabular}{c}}
\newcommand{\etab}{\end{tabular}}

\newtheorem{property}{Property}
\newtheorem{definition}{Definition}

\usepackage{import}

\title{An analytic comparison of regularization methods for Gaussian Processes}
\date{}
\author[1, 2]{Hossein Mohammadi }
\author[2, 1]{Rodolphe Le Riche\thanks{Corresponding author: Ecole Nationale Supérieure des Mines de Saint Etienne, Institut H. Fayol, 158, Cours Fauriel, 42023 Saint-Etienne cedex 2 - France \\ Tel : +33477420023 \\ Email: leriche@emse.fr}} 
\author[1, 2]{Nicolas Durrande} 
\author[1, 2]{Eric Touboul} 
\author[1, 2]{Xavier Bay} 
\affil[1]{Ecole des Mines de Saint Etienne, H. Fayol Institute}
\affil[2]{CNRS LIMOS, UMR 5168}

\begin{document}
\maketitle
\vskip .5cm

\section*{Abstract}
Gaussian Processes (GPs) are a popular approach to predict the output of a parameterized experiment.
They have many applications in the field of Computer Experiments, in particular to perform 
sensitivity analysis, adaptive design of experiments and global optimization.
Nearly all of the applications of GPs require the inversion of a covariance matrix
that, in practice, is often ill-conditioned. Regularization methodologies are then employed with 
consequences on the GPs that need to be better understood. 

The two principal methods to deal with ill-conditioned covariance matrices are 
\textit{i)} pseudoinverse and 
\textit{ii)} adding a positive constant to the diagonal (the so-called nugget regularization). 
The first part of this paper provides an algebraic comparison of PI and nugget regularizations.
Redundant points, responsible for covariance matrix singularity, are defined. 
It is proven that
pseudoinverse regularization, contrarily to nugget regularization, averages the output values and makes the variance zero at redundant points.
However, pseudoinverse and nugget regularizations become equivalent as the nugget value vanishes.
A measure for data-model discrepancy is proposed which serves for choosing a regularization technique.

In the second part of the paper, a distribution-wise GP is introduced that interpolates 
Gaussian distributions instead of data points. Distribution-wise GP can be seen as
an improved regularization method for GPs.

\vskip\baselineskip
\noindent Keywords: ill-conditioned covariance matrix; Gaussian Process regression; Kriging; Regularization. \\

\subsection*{Nomenclature} 

\subsubsection*{Abbreviations} 

CV, Cross-validation \\ \noindent 
discr, model-data discrepancy \\ \noindent 
GP, Gaussian Process \\ \noindent 
ML, Maximum Likelihood \\ \noindent 
PI, Pseudoinverse 

\subsubsection*{Greek symbols} 

$\tau^2$, nugget value. \\ \noindent 
$\Delta$, the difference between two likelihood functions. \\ \noindent 
$\kappa$, condition number of a matrix. \\ \noindent 
$\kappa_{max}$, maximum condition number after regularization. \\ \noindent 
$\lambda_i$, the $i$th largest eigenvalue of the covariance matrix. \\ \noindent 
$\mu(.)$, Gaussian process mean. \\ \noindent 
$\sigma^2$, process variance. \\ \noindent 
{\boldmath$\Sigma$}, diagonal matrix made of covariance matrix eigenvalues.  \\ \noindent 
$\eta$, tolerance of pseudoinverse. \\ \noindent 
$\theta_i$, characteristic length-scale in dimension $i$. 

\subsubsection*{Latin symbols} 
$\textbf{c}$, vector of covariances between a new point and the design points $\textbf X$. \\ \noindent 
$\textbf{C}$, covariance matrix. \\ \noindent 
$\textbf{C}^i$, $i$th column of $\textbf{C}$. \\ \noindent 
$\textbf{e}_i$, $i$th unit vector. \\ \noindent 
$f:\mathbb{R}^d \rightarrow \mathbb{R}$, true function, to be predicted.  \\ \noindent 
$\textbf{I}$, identity matrix. \\ \noindent
$K$, kernel or covariance function.  \\ \noindent 
$m(.)$, kriging mean. \\ \noindent 
$n$, number of design points. \\ \noindent 
$N$, number of redundant points. \\ \noindent 
$\textbf{P}_{Im}$, orthogonal projection matrix onto the image space of a matrix (typically $\textbf{C}$). \\ \noindent 
$\textbf{P}_{Nul}$, orthogonal projection matrix onto the null space of a matrix (typically $\textbf{C}$). \\ \noindent
$\textbf{R}$, correlation matrix. \\ \noindent  
$r$, rank of the matrix $\textbf{C}$. \\ \noindent 
$s^2(.)$, kriging variance. \\ \noindent
$s_i^2$, variance of response values at $i$-th repeated point. \\ \noindent 
$\textbf{V}$, column matrix of eigenvectors of $\textbf{C}$ associated to strictly positive eigenvalues. \\ \noindent 
$\textbf{W}$, column matrix of eigenvectors of $\textbf{C}$ associated to zero eigenvalues. \\ \noindent 
$\textbf{X}$, matrix of design points. \\ \noindent 
$Y(.)$, Gaussian process. \\ \noindent 
$\textbf{y}$, vector of response or output values. \\ \noindent 
$\overline{y_i}$, mean of response values at $i$-th repeated point.

\section{Introduction}
Conditional Gaussian Processes, also known as kriging models, are commonly used for predicting from a set of spatial observations. 
Kriging performs a linear combination of the observed response values. 
The weights in the combination depend only, through a covariance function, on the locations of the points where one wants to predict and the locations of the observed points \cite{cressie1993, sacks1989, welch1992, GPML}.
%
%
We assume in this work that the location of the observed points and the covariance 
function are given a priori, which is the default situation when using algorithms performing
adaptive design of experiments \cite{sequential_DoE}, global sensitivity analysis \cite{sensitivity_analysis} and global optimization \cite{jones1998}. 

Kriging models require the inversion of a covariance matrix which is made of the covariance function evaluated at every pair of observed locations. 
In practice, anyone who has used a kriging model has experienced one of the circumstances when the covariance matrix is ill-conditioned, hence not reliably invertible by a numerical method. 
This happens when observed points are too close to each other, or more generally when the covariance function makes the information provided by observations redundant.

\vskip\baselineskip

In the literature, various strategies have been employed to avoid such degeneracy of the covariance matrix. 
A first set of approaches proceed by controlling the locations of design points (the Design of Experiments or DoE).
The influence of the DoE on the condition number of the covariance matrix has been investigated in \cite{factorial}. 
\cite{rennen} proposes to build kriging models from a uniform subset of design points to improve the condition number. 
In \cite{osborne2009}, new points are taken suitably far from all existing data points to guarantee a good conditioning.

Other strategies select the covariance function so that the covariance matrix remains well-conditioned.
In \cite{sixfactor} for example, the influence of all kriging parameters  on the condition number, including the covariance function, is discussed. 
Ill-conditioning also happens in the related field of linear regression with the Gauss-Markov matrix $\bf{\Phi}^\top \bf{\Phi}$ that needs to be inverted, where $\bf{\Phi}$ is the matrix of basis functions evaluated at the DoE. In regression, work has been done on diagnosing ill-conditioning and the solution typically involves working on the definition of the basis functions to recover invertibility \cite{belsley1991}. 
The link between the choice of the basis functions and the choice of the covariance functions is given by Mercer's theorem, \cite{GPML}.

Instead of directly inverting the covariance matrix, an iterative method has been proposed in \cite{GibbstPhd} to solve the kriging equations and 
avoid numerical instabilities. 

\vskip\baselineskip

The two standard solutions to overcome the ill-conditioning of covariance matrices are the pseudoinverse (PI) and the ``nugget" regularizations. 
They have a wide range of applications because, contrarily to the methods mentioned above, they can be used a posteriori in computer experiments algorithms without 
major redesign of the methods. This is the reason why most kriging implementations contain PI or nugget regularization.

The singular value decomposition and the idea of pseudoinverse have already been suggested in \cite{jones1998} in relation with Gaussian Processes (GPs). 
The Model-Assisted Pattern Search (MAPS) software \cite{siefert} relies on an implementation of the pseudoinverse to invert the covariance matrices. 

The most common approach to deal with covariance matrix ill-conditioning
is to introduce a ``nugget'' \cite{booker1998b, santer2003, neal1997, andrianakis2012}, that is to say add a small positive scalar to the diagonal. As a matrix regularization method, it is also known as Tikhonov regularization or Ridge regression.
The popularity of the nugget regularization may be due to its simplicity and to its interpretation as the variance of a noise on the observations.
The value of the nugget term can be estimated by maximum likelihood (ML). 
It is reported in \cite{pepelyshev2010} that the presence of a nugget term significantly changes the modes of the likelihood function of a GP. 
Similarly in \cite{gramacy2009}, the authors have advocated a nonzero nugget term in the design and analysis of their computer experiments.
They have also stated that estimating a nonzero nugget value may improve some statistical properties 
of the kriging models such as their predictive accuracy \cite{gramacy2012}. 
In contrast, some references like \cite{ranjan2011} recommend that the magnitude of nugget remains as small as possible to preserve the interpolation property.

\vskip\baselineskip

Because of the diversity of arguments regarding GP regularization, we feel that there is a need to provide 
analytical explanations on the effects of the main approaches. 
The paper starts by detailing how covariance matrices associated to GPs become singular, which leads to 
the definition of redundant points. Then, new results are provided regarding the analysis and comparison of pseudoinverse and nugget kriging regularizations.
The analysis is made possible by approximating ill-conditioned covariance matrices with the neighboring truly singular covariance matrices. 
The paper finishes with the description of a new regularization method associated to distribution-wise GPs.

\section{Kriging models and degeneracy of the covariance matrix} \label{kriging_section}
\subsection{Context: conditional Gaussian processes}
This section contains a summary of conditional GP concepts and notations. Readers who are familiar with GP may want to proceed to the next section (\ref{sec-degeneracy}).

Let $f$ be a real-valued function defined over $D\subseteq \mathbb{R}^d$. Assume that the values of $f$ are known at a limited set of points 
called design points. One wants to infer the value of this function elsewhere. Conditional GP is one of the most important technique for this purpose \cite{sacks1989, welch1992}. 

A GP defines a distribution over functions. Formally, a GP indexed by $D$ is a collection of random variables $\left(Y(\textbf{x}); \textbf{x}\in D \right)$ such that for any $n \in \mathbb{N}$ and any $\textbf{x}^1, ..., \textbf{x}^n \in D$, $\left(Y(\textbf{x}^1), ..., Y(\textbf{x}^n) \right)$ follows a multivariate Gaussian distribution. 
The distribution of the GP is fully characterized by a mean function $\mu(\textbf{x}) = \mathbb{E}(Y(\textbf{x}))$
and a covariance function $K(\textbf{x},\textbf{x'}) = \Cov(Y(\textbf{x}),Y(\textbf{x'}))$ \cite{GPML}.

The choice of kernel plays a key role in the obtained kriging model. In practice, a parametric family of kernels is selected (e.g., Mat\'ern, polynomial, exponential) and then the unknown kernel parameters are estimated from the observed values. For example, a separable squared exponential kernel is expressed as
\begin{eqnarray} \label{E:gauss_kernel}
K\left(\textbf{x},\textbf{x}^\prime\right)=\sigma^2 \prod\limits_{i=1}^d \exp \left(-\frac{\mid x_i- x^\prime_i\mid ^2}{2\theta^2_i} \right).
\end{eqnarray}
\noindent In the above equation, $\sigma^2$ is a scaling parameter known as process variance and $x_i$ is the $i$th component of $\textbf{x}$. 
The parameter $\theta_i$ is called length-scale and determines the correlation length along coordinate $i$.
It should be noted that $\Cov\left(Y(\textbf{x}),Y(\textbf{x}^\prime)\right)$ in Equation~(\ref{E:gauss_kernel}) is only a function of the difference between $\textbf{x}$ and $\textbf{x}^\prime$. A GP with this property is said to be stationary, otherwise it is nonstationary. Interested readers are referred to \cite{GPML} for further information about GPs and kernels.

Let $Y(\textbf{x})_{\textbf{x} \in D}$ be a GP with kernel $K(.,.)$ and zero mean ($\mu(.)=0$). 
$\textbf{X}=\left(\textbf{x}^1, ..., \textbf{x}^n \right)$ denotes the $n$ data points where the samples are taken and the corresponding response values are $\textbf{y}=\left(y_1, ..., y_n \right)^\top = \left(f(\textbf{x}^1), ..., f(\textbf{x}^n) \right)^\top$. 
The posterior distribution of the GP $\left(Y(\textbf{x})\right)$ knowing it interpolates the data points is still Gaussian with mean and covariance \cite{GPML}
\begin{eqnarray}
m(\textbf{x}) &=& \mathbb{E}(Y(\textbf{x}) | Y(\textbf{X})=\textbf{y}) ~=~\textbf{c}(\textbf{x})^\top \textbf{C}^{-1}\textbf{y}~, \label{E:kriging_mean}\\
c(\textbf{x},\textbf{x}^\prime) &=& \Cov(Y(\textbf{x}),Y(\textbf{x}^\prime) | Y(\textbf{X})=\textbf{y})  \nonumber\\
~ &=& K(\textbf{x},\textbf{x}^\prime) - \textbf{c}(\textbf{x})^\top \textbf{C}^{-1}\textbf{c}(\textbf{x}^\prime) \label{E:kriging_variance}~,
\end{eqnarray}
where $\textbf{c}(\textbf{x})=\left(K(\textbf{x},\textbf{x}^1), ..., K(\textbf{x},\textbf{x}^n)\right)^\top$ is the vector of covariances between a new point $\textbf{x}$ and the $n$ already observed sample points. 
The $n\times n$ matrix $\textbf{C}$ is a covariance matrix between the data points and its elements are defined as $\textbf{C}_{i,j}=K \left(\textbf{x}^i,\textbf{x}^j \right)=\sigma^2 \textbf{R}_{i,j}$, where $\textbf{R}$ is the correlation matrix. Hereinafter, we call $m(\textbf{x})$ and $v(\textbf{x}) = c(\textbf{x},\textbf{x})$ the kriging mean and variance, respectively.

One essential question is how to estimate the unknown parameters in the covariance function. Typically, the values of the model parameters (i.e., $\sigma$ and the $\theta_i$'s) are learned via maximization of the likelihood. The likelihood function of the unknown parameters given observations $\textbf{y}=\left(y_1, ..., y_n \right)^\top$ is defined as follows:
\begin{eqnarray}
L(\textbf{y}\vert \bm{\theta}, \sigma^2)=\frac{1}{(2\pi)^{n/2} \vert \textbf{C} \vert ^{1/2}} \exp\left(-\frac{\textbf{y}^\top \textbf{C}^{-1}\textbf{y}}{2} \right).
\end{eqnarray}
\noindent In the above equation, $\vert \textbf{C} \vert$ indicates the determinant of the covariance matrix $\textbf{C}$ and {\boldmath$\theta$}$=(\theta_1, ..., \theta_d)^\top$ is a vector made of the length-scales in each dimension.
It is usually more convenient to work with the natural logarithm of the likelihood function that is:
\begin{eqnarray} \label{E:loglik}
\ln L(\textbf{y}\vert \bm{\theta}, \sigma^2)=-\frac{n}{2}\ln(2\pi) - \frac{1}{2}\ln\vert \textbf{C} \vert - \frac{1}{2}\textbf{y}^\top \textbf{C}^{-1}\textbf{y}.
\end{eqnarray}
The ML estimator of the process variance $\sigma^2$ is
\begin{eqnarray}
\hat{\sigma}^2= \frac{1}{n} \textbf{y}^\top \tilde{\textbf{K}}^{-1}\textbf{y},
\end{eqnarray}  
and if it is inserted in (\ref{E:loglik}), it yields (minus) the concentrated log-likelihood,
\begin{eqnarray} \label{E:-2loglik}
-2\ln L(\textbf{y}\vert \bm{\theta}, \sigma^2)= n\ln(2\pi) + n\ln\hat{\sigma}^2 + \ln\vert \tilde{\textbf{K}}\vert + n.
\end{eqnarray} 
Finally, {\boldmath$\theta$} is estimated by numerically minimizing Equation~(\ref{E:-2loglik}).

\subsection{Degeneracy of the covariance matrix} \label{degeneracy_cov.mat}
\label{sec-degeneracy}

Computing the kriging mean (Equation~(\ref{E:kriging_mean})) or (co)variance (Equation~(\ref{E:kriging_variance})) or even samples of GP trajectories, 
requires inverting the covariance matrix $\textbf{C}$. 
In practice, the covariance matrix should not only be invertible, but also well-conditioned.  A matrix is said to be near singular or ill-conditioned or degenerated if its condition number is too large. For covariance matrices, which are symmetric and positive semidefinite, the condition number $\kappa(\textbf{C})$ is the ratio of the largest to the smallest eigenvalue. Here, we assume that $\kappa(\textbf{C}) \rightarrow \infty$ is possible. 

\vskip\baselineskip

There are many situations where the covariance matrix is near singular. The most frequent and easy to understand case is 
when some data points are too close to each other, where closeness is measured with respect to the metric induced by the covariance function. 
This is a recurring issue in sequential DoEs like the EGO algorithm \cite{jones1998} where the search points tend to pile up around the points of interest such as the global optimum \cite{ranjan2011}. When this happens, the resulting covariance matrix is no longer numerically invertible because some columns are almost identical. 

\vskip\baselineskip

Here, to analyze PI and nugget regularizations, 
we are going to consider matrix degeneracy pushed to its limit, that is $\textbf{C}$ is 
mathematically non-invertible (i.e., it is rank deficient).
Non invertibility happens if a linear dependency exists between columns (or rows) of $\textbf{C}$.
Section~\ref{app-redundant} provides a collection of examples where the covariance matrix is not invertible with calculation details that will become clear later.
Again, the easiest to understand and the most frequent occurrence of $\textbf{C}$'s rank deficiency is when some of the data 
points tend towards each other until they are at the same $\textbf x^i$ position. They form \emph{repeated} points, 
the simplest example of what we more generally call redundant points which will be formally defined shortly.
Figure~\ref{fig-repeated} in Section~\ref{app-redundant} is an example of repeated points.
Repeated points lead to strict non-invertibility of $\textbf{C}$ since the corresponding columns are identical.
The special case of repeated points will be instrumental 
in understanding some aspects of kriging regularization in Sections \ref{PI_averaging} and \ref{nugget_and_ML} 
because the eigenvectors of the covariance matrix associated to eigenvalues equal to zero are known. 

\vskip\baselineskip

The covariance matrix of GPs may loose invertibility even though the data points are not close to each other in 
Euclidean distance. This occurs for example with additive GPs for which the kernel is the sum of kernels defined in each dimension, 
$K(\textbf{x}, \textbf{x}')=\sum\limits_{i=1}^d K_i(x_i, x_i')$.
The additivity of a kernel may lead to linear dependency in some columns of the covariance matrix. 
For example, in the DoE shown in Figure \ref{additive_kernel}, only three of the first four points which form a rectangle provide independent information
in the sense that the GP response at any of the four points in fully defined by the response at the three other points.
This is explained by a linear dependency between the first four columns,
$\textbf{C}^4=\textbf{C}^3+\textbf{C}^2-\textbf{C}^1$, which comes from the additivity of the kernel and the rectangular design \cite{durrande2012}:
\begin{equation*} \label{E:add_kernel0}
\textbf{C}^4_i ~=~\Cov(x^i_1,x^4_1) + \Cov(x^i_2,x^4_2) ~=~ \Cov(x^i_1,x^2_1) + \Cov(x^i_2,x^3_2)~,
\end{equation*} 
and completing the covariances while accounting for $x^2_2=x^1_2$, $x^3_1=x^1_1$, yields
\begin{equation*} \label{E:add_kernel}
\textbf{C}^4_i~=~\Cov(\textbf{x}^i, \textbf{x}^3)+\Cov(\textbf{x}^i, \textbf{x}^2)-\Cov(\textbf{x}^i, \textbf{x}^1) ~=~ \textbf{C}^3_i + \textbf{C}^2_i - \textbf{C}^1_i~.
\end{equation*}
Note that if the measured outputs $y^1$, \ldots, $y^4$ are not additive ($y^4 \ne y^2 + y^3 - y^1$), none of the four 
measurements can be easily deleted without loss of information, 
hence the need for the general regularization methods that will be discussed later. \\
Periodic kernels may also yield non-invertible covariance matrices although data points are far from each other.
This is illustrated in Figure~\ref{fig-periodic}  
where points 1 and 2, and points 3 and 4, provide the same information as they are one period away from each other. 
Thus, $\textbf{C}^1=\textbf{C}^2$ and $\textbf{C}^3=\textbf{C}^4$. \\
Our last example comes from the dot product (or linear) kernel (cf. Section~\ref{sec-linear}). Because the GP trajectories and mean are linear, no uncertainty is left in the model when the number of data points $n$ reaches $d+1$ and when $n>d+1$ the covariance matrix is no longer invertible.

\subsection{Eigen analysis and definition of redundant points} \label{sec-redundant}
We start by introducing our notations for the eigendecomposition of the covariance matrix.
Let the $n \times n$ covariance matrix $\textbf{C}$ have rank $r$, $r \le n$. 
A covariance matrix is positive semidefinite, thus its eigenvalues are greater than or equal to zero. 
The eigenvectors associated to strictly positive eigenvalues are denoted $\textbf{V}^i$, $i=1,\ldots,r$, 
and those associated to null eigenvalues are $\textbf{W}^i$, $i=1,\ldots,(n-r)$, that is 
$\textbf{C} \textbf{V}^i = \lambda_i \textbf{V}^i$ where $\lambda_i>0$ and $\textbf{C} \textbf{W}^i =\textbf{0}$.
The eigenvectors are grouped columnwise into the matrices $\textbf{V} = [\textbf{V}^1 , \ldots , \textbf{V}^r ]$ and
$\textbf{W} = [\textbf{W}^1 , \ldots , \textbf{W}^{n-r} ]$. 
In short, the eigenvalue decomposition of the covariance matrix $\textbf{C}$ obeys	
\begin{align}\label{E:EVD}
\textbf{C} ~=~ [\textbf{V}~\textbf{W}] \, \mathbf{\Sigma} \, [\textbf{V}~\textbf{W}]^\top,
\end{align}
\noindent where {\boldmath$\Sigma$} is a diagonal matrix containing the eigenvalues of $\textbf{C}$, $\lambda_1 \ge \lambda_2 \ge \ldots \ge \lambda_r > 0$ and $\lambda_{r+1}=\ldots=\lambda_n=0$. 
$\textbf{V}$ spans the image space and $\textbf{W}$ spans the null space of $\textbf{C}$, $Im(\textbf{C})$ and $Null(\textbf{C})$, respectively. 
$[\textbf{V}~\textbf{W}]$ is an orthogonal matrix, 
\begin{equation}
[\textbf{V}~\textbf{W}]^\top [\textbf{V}~\textbf{W}]= [\textbf{V}~\textbf{W}] [\textbf{V}~\textbf{W}]^\top 
=\textbf{V} \textbf{V}^\top + \textbf{W} \textbf{W}^\top = \textbf{I}~. 
\label{eq-vwid}
\end{equation}
$\textbf{V} \textbf{V}^\top$ is the orthogonal projection matrix onto $Im(\textbf{C})$. 
Similarly, $\textbf{W} \textbf{W}^\top$ is the orthogonal projection matrix onto $Null(\textbf{C})$.
For a given matrix $\textbf{C}$, the eigenvectors $\textbf{W}^i$ are not uniquely defined 
because any linear combination of them is also an eigenvector associated to a null eigenvalue. 
However, the orthogonal projection matrices onto the image and null spaces of $\textbf{C}$ are unique and 
will be cornerstones in the definition of redundant points.

\vskip\baselineskip

Before formally defining redundant points, we present the examples of singular covariance matrices 
of  Section~\ref{app-redundant}. These examples are 
two dimensional to allow for a graphical representation. The kernels, designs of points, 
eigenvalues and eigenvectors and 
the $\textbf{V} \textbf{V}^\top$ projection matrix are given. \\
The first example detailed in Section \ref{app-redundant} has two groups of repeated data points (points 1, 2 and 6, on the one hand, points 3 and 4, on the other hand),
in which there are 3 redundant, points. The covariance matrix has 3 null eigenvalues. 
It should be noted that the off-diagonal coefficients of 
the $\textbf{V} \textbf{V}^\top$ projection matrix associated to the indices of repeated points are not 0. \\
Figure \ref{fig-add1} shows how additive kernels may generate singular covariance matrices:
points 1, 2, 3 and 4 are arranged in a rectangular pattern which makes columns 1 to 4 linearly dependent 
(as already explained in Section~\ref{sec-degeneracy}). The additive property makes any one of the 4 points 
of a rectangular pattern redundant in that the value of the GP there is uniquely set by the 
knowledge of the GP at the 3 other points.
The same stands for points 5 to 8. Two points are redundant (1 in each rectangle) and there are two null eigenvalues.
Again, remark how the off-diagonal coefficients of $\textbf{V} \textbf{V}^\top$ associated to the points of the rectangles are not zero.
Another example of additivity and singularity is depicted in Figure~\ref{fig-add2}: 
although the design points are not set in a rectangular pattern, 
there is a shared missing vertex between two orthogonal triangles so that, because of additivity, the value at this missing 
vertex is defined twice. In this case, there is one redundant point, one null eigenvalue, and all the points of the design are
coupled: all off-diagonal terms in $\textbf{V} \textbf{V}^\top$ are not zero.\\
Finally, Figure~\ref{fig-periodic} is a case with a periodic kernel and a periodic pattern of points so that 
points 1 and 2 provide the same information, and similarly with points 3 and 4. 
There are 2 null eigenvalues, and the (1,2) and (3,4) off-diagonal terms in $\textbf{V} \textbf{V}^\top$ are not zero. \\

In general, we call \textit{redundant} the set of data points that make the covariance 
matrix non-invertible by providing linearly dependent information.
\begin{definition}[Redundant points set]
\mbox{~}\\
Let {\normalfont $\textbf{C}$} be a {\normalfont $n \times n$} positive semidefinite covariance matrix of rank {\normalfont $r$} whose 
generic term {\normalfont $\textbf{C}_{i,j}$} is associated to the points {\normalfont $\textbf{x}^i$} and {\normalfont $\textbf{x}^j$} through
{\normalfont $\textbf{C}_{i,j} = \Cov(Y(\textbf{x}^i),Y(\textbf{x}^j))$}.
{\normalfont $\textbf{V}$} is the {\normalfont $n \times r$} matrix of the eigenvectors associated to strictly positive eigenvalues.
Two points $i$ and $j$, $i\ne j$, are called redundant when {\normalfont $(\textbf{V} \textbf{V}^\top)_{ij} \ne 0$}.
The set of redundant points, {\normalfont $R$}, is made of all point indices that are redundant.
\end{definition}
Redundant points could be equivalently defined with the $\textbf{W}$ matrix since, 
from Equation~(\ref{eq-vwid}), $\textbf{V} \textbf{V}^\top$ and $\textbf{W} \textbf{W}^\top$
have the same non-zero off-diagonal terms with opposite signs. Subsets of redundant points are also redundant.
The \textit{degree of redundancy} of a set of points $R$ is the number of zero eigenvalues 
of the covariance matrix restricted to the points in $R$, i.e., 
$[\textbf{C}_{ij}$] for all $(i,j) \in R^2$. The degree of redundancy is the number of points 
that should be removed from $R$ to recover invertibility of the covariance restricted 
to the points in $R$. When $r=n$, $\textbf{C}$ is invertible and there is no redundant point.
An interpretation of redundant points will be made in the next Section on pseudoinverse regularization.

In the repeated points example of Section~\ref{app-redundant}, the two largest 
redundant points sets are $\{1,2,6\}$ and $\{3,4\}$ with degrees of redundancy 2 and 1, respectively.
The first additive example has two sets of redundant points,
$\{1,2,3,4\}$ and $\{5,6,7,8\}$ each with a degree of redundancy equal to 1. 
In the second additive example, all the points are redundant with 
a degree equal to 1. 
In the same section, the periodic case has two sets of redundant points of degree 1, $\{1,2\}$ and $\{3,4\}$.
With the linear kernel all data points are redundant and in the given example where $n=d+2$ the 
degree of redundancy is 1.

\section{Pseudoinverse regularization}\label{sec:pi}
\subsection{Definition} \label{kriging_by_PI}
In this Section, we state well-known properties of pseudoinverse matrices without proofs (which can be found, e.g., in \cite{benisrael}) and apply them to the kriging Equations (\ref{E:kriging_mean}) and (\ref{E:kriging_variance}). 
Pseudoinverse matrices are generalizations of the inverse matrix. 
The most popular pseudoinverse is the \textit{Moore–Penrose pseudoinverse} which is hereinafter referred to as pseudoinverse. 

When $\textbf{C}^{-1}$ exists (i.e., $\textbf{C}$ has full rank, $r=n$), we denote as {\boldmath$\beta$} the term $\textbf{C}^{-1}\textbf{y}$ of the kriging mean formula, Equation~(\ref{E:kriging_mean}). 
More generally, when $\textbf{C}$ is not a full rank matrix, we are interested in the vector {\boldmath$\beta$} that simultaneously minimizes\footnote{Indeed, in this case the minimizer of $\Vert \textbf{C}\bm{\beta}-\textbf{y}\Vert_2$ is not unique but defined up to any sum with a vector in the $Null(\textbf{C})$} $\Vert \textbf{C}\bm{\beta}-\textbf{y}\Vert_2$ and $\Vert\bm{\beta}\Vert_2$. This solution is unique and obtained by {\boldmath$\beta$}$^{PI}= \textbf{C} ^\dagger \textbf{y}$ where $\textbf{C}^\dagger$ is the pseudoinverse of $\textbf{C}$.
Each vector {\boldmath$\beta$} can be uniquely decomposed into
\begin{eqnarray}
\bm{\beta}=\bm{\beta}^{PI} + \bm{\beta}_{Null(\textbf{C})},
\end{eqnarray}
\noindent where {\boldmath$\beta$}$^{PI}$ and {\boldmath$\beta$}$_{Null(\textbf{C})}$ belong to the image space and the null space of the covariance matrix, respectively. The decomposition is unique since, $\textbf{C}$ being symmetric, $Im(\textbf{C})$ and $Null(\textbf{C})$ have no intersection.

\begin{figure}[htpb]
\begin{center}
\includegraphics[width=\textwidth]{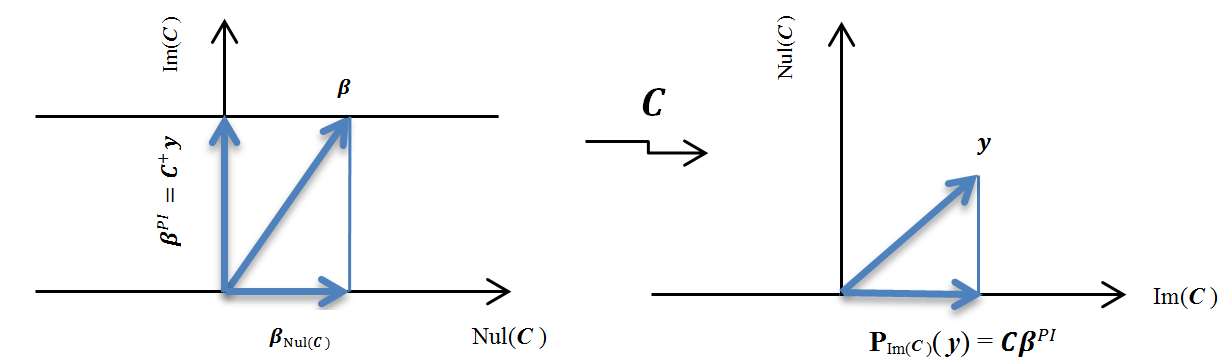}
 \caption{
Geometrical interpretation of the Moore-Penrose pseudoinverse. 
In the left picture, infinitely many vectors {\boldmath$\beta$} are solutions to the system $\textbf{C}${\boldmath$\beta$}$=\textbf{y}$. 
But the minimum norm solution is $\textbf{C}^\dagger \textbf{y}$. The right picture shows the orthogonal projection of $\textbf{y}$ onto the image space of $\textbf{C}$, $\textbf{P}_{Im(\textbf{C})}(\textbf{y})$, which is equal to $\textbf{CC}^\dagger \textbf{y}$ (Property \ref{prop-PIproj}).
\label{PI_geometry}
}
\end{center}
\end{figure}

The pseudoinverse of $\textbf{C}$ is expressed as
\begin{eqnarray} \label{E:PI}
\textbf{C}^\dagger= \left[ \textbf{V}~ \textbf{W} \right] \begin{bmatrix}
diag\left(\frac{1}{\bm{\lambda}}\right)_{r\times r} & \textbf{0}_{r \times (n-r)} \\[0.3em]
\textbf{0}_{(n-r)\times r} & \textbf{0}_{(n-r)\times (n-r)} 
\end{bmatrix} 
\left[ \textbf{V} ~ \textbf{W} \right]^\top,
\end{eqnarray}
where $diag(\frac{1}{\bm{\lambda}})$ is a diagonal matrix with $\frac{1}{\lambda_i}$, $i=1,\ldots,r$, as diagonal elements. So {\boldmath$\beta$}$^{PI}$ reads
\begin{eqnarray} \label{E:beta_PI}
\bm{\beta}^{PI}= \sum \limits_{i=1}^r \frac{\left(\textbf{V}^i \right)^\top \textbf{y}}{\lambda_i} \textbf{V}^i.
\end{eqnarray}
Equation~(\ref{E:beta_PI}) indicates that {\boldmath$\beta$}$^{PI}$ is in the image space of $\textbf{C}$, because it is a linear combination of eigenvectors associated to positive eigenvalues. A geometrical interpretation of {\boldmath$\beta$}$^{PI}$ and pseudo-inverse is given in Figure~\ref{PI_geometry}.

From now on, ``PI kriging'' will be a shorthand for ``kriging with pseudo-inverse regularization''.
The PI kriging mean (Equation~(\ref{E:kriging_mean})) can be written as
\begin{eqnarray} \label{E:kriging_mean_PI}
m^{PI}(\textbf{x})=\textbf{c}(\textbf{x})^\top \sum \limits_{i=1}^r \frac{\left(\textbf{V}^i \right)^\top \textbf{y}}{\lambda_i} \textbf{V}^i.
\end{eqnarray} 
Similarly, the kriging covariance (\ref{E:kriging_variance}) regularized by PI is,
\begin{equation} \label{E:kriging_variance_PI}
\begin{split}
c^{PI}(\textbf{x},\textbf{x}^\prime) &= K(\textbf{x},\textbf{x}^\prime) - \textbf{c}(\textbf{x})^\top\sum \limits_{i=1}^r \left( \frac{\left(\textbf{V}^i \right)^\top \textbf{c}(\textbf{x}^\prime)} {\lambda_i}\right)\textbf{V}^i \\
& = K(\textbf{x},\textbf{x}^\prime) - \sum \limits_{i=1}^r \frac{\left( \left(\textbf{V}^i \right)^\top \textbf{c}(\textbf{x})\right) \left( \left(\textbf{V}^i \right)^\top \textbf{c}(\textbf{x}^\prime)\right) }{\lambda_i}.
\end{split}
\end{equation}

\subsection{Properties of kriging regularized by PI} \label{PI_averaging}
\begin{figure}[htpb]
\centering
\includegraphics[width=0.5\textwidth]{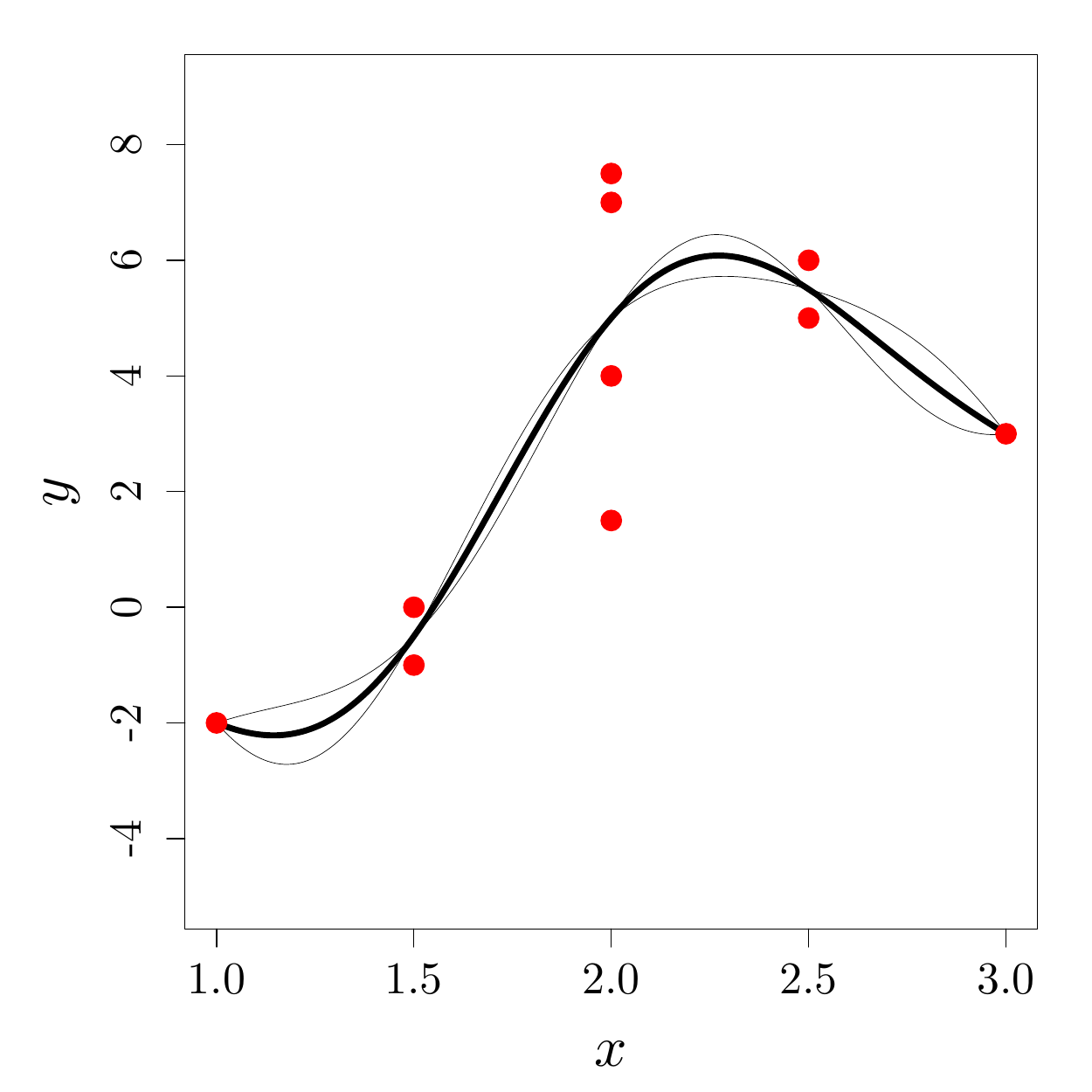}
\caption{
Kriging mean $m^{PI}(x)$ (thick line) and prediction intervals $m^{PI}(x) \pm 2\sqrt{v^{PI}(x)}$ (thin lines).
Kriging mean using pseudoinverse goes exactly through the average of the outputs. 
The observed values are $\textbf{y}=(-2, -1, 0, 1.5, 4, 7, 7.5, 6, 5, 3)^\top$. $m^{PI}(1.5)=-0.5$, $m^{PI}(2)=5$, and $m^{PI}(2.5)=5.5$.
Note that $v^{PI}$ is zero at redundant points.
\label{PIaveraging}
}
\end{figure}

The kriging mean with PI averages the outputs at redundant points. 
Before proving this property, let us illustrate it with the simple example
of Figure \ref{PIaveraging}: there are redundant points at $\textbf{x}=1.5$, $\textbf{x}=2$ and $\textbf{x}=2.5$.
We observe that the kriging mean with PI 
regularization is equal to the mean of the outputs, $m^{PI}(1.5) = -0.5 = (-1+0)/2$, $m^{PI}(2) = 5 = (1.5+4+7+7.5)/4$ and $m^{PI}(2.5) = 5.5 = (5+6)/2$.
The PI averaging property is due to the more abstract fact that PI projects the observed $\textbf y$ onto the image space of $\textbf C$.

\begin{property}[PI as projection of outputs onto {\normalfont $Im(\textbf{C})$}] \label{prop-PIproj}
\mbox{~}\\
The kriging prediction with PI regularization at {\normalfont $\textbf X$} is the projection of the observed outputs onto the image space of the covariance matrix, {\normalfont $Im(\textbf{C})$}. 
\end{property}
\noindent \textit{Proof}: The kriging means at all design points is given by
\begin{equation}
m^{PI}(\textbf{X})=\textbf{CC}^\dagger \textbf{y}~.
\label{eq-mpi}
\end{equation}
Performing the eigendecompositions of the matrices, one gets,
\begin{eqnarray}
m^{PI}(\textbf{X}) &=& \left[ \textbf{V} ~ \textbf{W} \right] 
\left[ \begin{array}{cc} diag(\bm{\lambda}) & \textbf{0} \\ \textbf{0} & \textbf{0} \\ \end{array}\right]
\left[ \begin{array}{c} \textbf{V}^\top \\ \textbf{W}^\top \end{array}\right]
\left[ \textbf{V} ~ \textbf{W} \right] 
\left[ \begin{array}{cc} diag(\frac{1}{\bm{\lambda}}) & \textbf{0} \\ \textbf{0} & \textbf{0} \\ \end{array}\right]
\left[ \begin{array}{c} \textbf{V}^\top \\ \textbf{W}^\top \end{array}\right] \textbf{y} 
 \nonumber \\
~ &=&  \textbf{V} \textbf{V}^\top \textbf{y} \label{eq-PIVVT} \end{eqnarray}
The matrix 
\begin{equation}
\label{eq-Pim}
\textbf{P}_{Im(\textbf C)} ~=~ \textbf{V} \textbf{V}^\top  ~=~ (\textbf{I} - \textbf{W} \textbf{W}^\top)
\end{equation}
is the orthogonal projection onto the image space of $\textbf{C}$ because it holds that
\begin{align*}
&\textbf{P}_{Im(\textbf C)}=\textbf{P}_{Im(\textbf C)}^\top ;  \\ 
&\textbf{P}_{Im(\textbf C)}^2=\textbf{P}_{Im(\textbf C)} ;\\ 
&\forall \textbf{v} \in  Im(\textbf{C}) ~,~ \textbf{P}_{Im(\textbf C)}\textbf{v}=\textbf{v} ; \\  
&\text{and } \forall \textbf{u} \in Null(\textbf{C})~,~ \textbf{P}_{Im(\textbf C)}\textbf{u}=\textbf{0} \quad \square 
\end{align*}

Redundant points can be further understood thanks to Property~\ref{prop-PIproj} and Equation~(\ref{eq-PIVVT}): 
points redundant with $\textbf x^i$ are points $\textbf x^j$ where the observations influences $m^{PI}(\textbf{x}^i)$.
The kriging predictions at the redundant data points $m^{PI}(\textbf{x}^i)$ and $m^{PI}(\textbf{x}^j)$ are not $\textbf{y}_i$ and $\textbf{y}_j$, as it happens at non-redundant points where the model is interpolating, but a linear combination of them.
The averaging performed by PI becomes more clearly visible in the important case of repeated points.

\begin{property}[PI Averaging Property for Repeated Points] \label{PI-averaging}
\mbox{~}\\
The PI kriging prediction at repeated points is the average of the outputs at those points.
\end{property}

\textit{Proof}: Suppose that there are $N$ repeated points at $k$ different locations with $N_i$ points at each repeated location, $\sum \limits_{i=1}^k N_i=N$, see Figure~\ref{F:redundant_points}.
The corresponding columns in the covariance matrix are identical,
\begin{equation*}
\textbf{C}=\left(\underbrace{\textbf{C}^1, ..., \textbf{C}^1}_{N_1 \text{ times} }, 
\ldots , \underbrace{\textbf{C}^k, ..., \textbf{C}^k}_{N_k \text{ times} }, \textbf{\textbf{C}}^{N+1}, ..., \textbf{\textbf{C}}^n \right).
\end{equation*}
In this case, the dimension of the image space, or rank of the covariance matrix, is $n - N + k$ and the dimension of the null space is equal to $\sum \limits_{i=1}^k (N_i - 1) = N - k$.

To prove this property we need to show that the matrix $\textbf{P}$ defined as
\begin{eqnarray}
\textbf{P} ~=~ \begin{pmatrix}
  \frac{\textbf{J}_{N_1}}{N_1} &  &  & & \\
    &  & \ddots & & 0 \\
  & 0 & & \frac{\textbf{J}_{N_k}}{N_k} & \\
   & & & & \textbf{I}_{n-N}
 \end{pmatrix},
\label{matrix_proj_image}
\end{eqnarray}
is the orthogonal projection matrix onto the image space of $\textbf{C}$, or $\textbf{P}=\textbf{P}_{Im(\textbf{C})}$. In Eq.~\ref{matrix_proj_image}, $\textbf{J}_{N_i}$ is the $N_i \times N_i$ matrix of ones and $\textbf{I}_{n-N}$ is the identity matrix of size $n-N$.
If $\textbf{P}=\textbf{P}_{Im(\textbf{C})}$, because of the unicity of the orthogonal projection and Property~\ref{prop-PIproj}, $m^{PI}(\textbf{X})$ is expressed as
\begin{eqnarray}\label{projection_IM}
m^{PI}(\textbf{X})=\textbf{P}_{Im(\textbf{C})}\textbf{y}=
\begin{bmatrix}
\overline{y_1} \\
\vdots \\
\overline{y_1} \\
\vdots \\
\overline{y_k} \\
\vdots \\
\overline{y_k} \\
y_{N+1} \\
\vdots \\
y_n
\end{bmatrix} ,
\end{eqnarray}
in which $\overline{y_i}= \frac{\sum \limits_{j=N_1+...+N_{i-1}+1}^{N_i}y_j}{N_i}$. It means that the PI kriging prediction at repeated points is the average of the outputs at those points.

It is easy to see that $\textbf{P}^\top=\textbf{P}$ and $\textbf{P}^2=\textbf{P}$. We now check the two remaining characteristic properties of projection matrices
\begin{enumerate}
        \item $\forall \textbf{u} \in Null(\textbf{C})~,~ \textbf{P}\textbf{u}=\textbf{0}$
        \item $\forall \textbf{v} \in  Im(\textbf{C}) ~,~ \textbf{P}\textbf{v}=\textbf{v}.$
\end{enumerate}
We first construct a set of non-orthogonal basis vectors of $Null(\textbf{C})$.
The basic idea is that when two columns of the covariance matrix $\textbf{C}$ are identical, e.g., the two first columns,
$\textbf{C} = \left( \textbf{C}^1 , \textbf{C}^1 , \ldots \right) $, then vector
$\textbf{u}^1 = (1, -1, 0,\ldots,0)^\top/\sqrt{2}$ belongs to $Null(\textbf{C})$ because
\begin{equation}
\textbf{C}^1 -\textbf{C}^1 = \textbf{C}\textbf{e}_1 - \textbf{C}\textbf{e}_2 = \textbf{C}(\underbrace{\textbf{e}_1 - \textbf{e}_2}_{\textbf{u}^1})=\textbf{0}.
\end{equation}
Generally, all such vectors can be written as
\begin{equation*}
\textbf{u}^j ~=~ \frac{\textbf{e}_{j+1} - \textbf{e}_{j}}{\sqrt{2}} ~,~ j = \sum_{l \le i-1} N^l +1, \ldots , \sum_{l \le i} N^l -1~,~ i=1,\ldots,k ~.
\end{equation*}
There are $ N-k = dim(Null(\textbf{C}))$ such $\textbf{u}^j$'s which are not orthogonal but linearly independent.
They make a basis of $Null(\textbf{C})$.
It can be seen that $\textbf{P}\textbf{u}^j=\textbf{0}~,~j=1, \dots, N-k$.
Since every vector in $Null(\textbf{C})$ is a linear combination of the $\textbf{u}^j$'s, the equation $\textbf{P}\textbf{u}=\textbf{0}$ holds for any vector in the null space of $\textbf{C}$ which proves the first characteristic property of the projection matrix.

The second property is also proved by constructing a set of vectors that span $Im(\textbf{C})$.
There are $n-N+k$ such vectors. The $k$ first vectors have the form
\begin{equation}
\textbf{v}^i = (\underbrace{0,\ldots, 0}_{N_1+\ldots+N_{i-1} \text{ times}} ,\quad \underbrace{1,\ldots,1}_{N_i \text{ times}}, \quad 0,\ldots,0)^\top/\sqrt{N_i}~,~ i=1, \dots, k.
\end{equation}
The $n-N$ other vectors are: $\textbf{v}^j=\textbf{e}_{j-k+N}~, j=k+1, \dots, n-N+k$.
Because these $n-N+k$ $~\textbf{v}^j$'s are linearly independent and perpendicular to the null space (to the above $\textbf{u}^j~,~ j=1, \dots, N-k$), they span $Im(\textbf{C})$.
Furthermore, $\textbf{P}\textbf{v}^j=\textbf{v}^j~,~ j=1, \dots, n-N+k$. 
The equation $\textbf{P}\textbf{v}=\textbf{v}$ is true for every $\textbf{v}\in Im(\textbf{C})$, therefore, $\textbf{P}$ is the projection matrix onto the image space of $\textbf{C}$ and the proof is complete.
$\square$

\begin{property}[Null variance of PI regularized models at data points] \label{PI-nulvar}
\mbox{~}\\
The variance of Gaussian processes regularized by pseudoinverse is zero at data points.
\end{property}
Therefore $v^{PI}(\cdot)$ is zero at redundant points. \\
\noindent \textit{Proof}: 
From Equation~(\ref{E:kriging_variance}), the PI kriging variances at all design points are 
\begin{equation*}
v^{PI}(\textbf{X})~=~c^{PI}(\textbf{X},\textbf{X})~=~K(\textbf{X},\textbf{X}) -  \textbf{c}(\textbf{X})^\top \textbf{C}^\dagger  \textbf{c}(\textbf{X}) ~=~ \textbf{C} - \textbf{C}^\top \textbf{C}^\dagger \textbf{C} = \textbf{C} - \textbf{C} ~=~0~,
\end{equation*} 
thanks to the pseudoinverse property \cite{strang}, $\textbf{C}\textbf{C}^\dagger \textbf{C}~=~\textbf{C}$. 
$\square$

\section{Nugget regularization } \label{section_nugget_regularization}

\subsection{Definition and covariance orthogonality property} \label{subsection_nugget_regularization}

When regularizing a covariance matrix by nugget, a positive value, $\tau^2$, is added to the main diagonal.
This corresponds to a probabilistic model with an additive white noise of variance $\tau^2$, ~$Y(\textbf x) ~\mid~ Y(\textbf x^i) + \varepsilon_i = \textbf y_i$, ~$i=1,\ldots,n$, ~where the $\varepsilon_i$'s are i.i.d. $\mathcal N(0,\tau^2) $.
Nugget regularization improves the condition number of the covariance matrix by increasing all the eigenvalues by $\tau^2$: 
if $\lambda_i$ is an eigenvalue of $\textbf{C}$, then $\lambda_i+\tau^2$ is an eigenvalue of $\textbf{C}+\tau^2 \textbf{I}$  and the eigenvectors remain the same (the proof is straightforward). The associated condition number is $\kappa(\textbf{C}+\tau^2 \textbf{I})=\frac{\lambda_{max}~+~\tau^2}{\lambda_{min}~+~\tau^2}$. 
The nugget parameter causes kriging to smoothen the data and become non-interpolating.

\begin{property}[Loss of interpolation in models regularized by nugget] \label{prop-lossinterp}
\mbox{~}\\
A conditional Gaussian process regularized by nugget has its mean no longer always equal 
to the output at data points, {\normalfont $m^{Nug}(\textbf{x}^i) ~\ne~y^i$, $i=1,n$}.
\end{property}
This property can be understood as follows. A conditional GP with invertible covariance matrix is interpolating because 
$c(\textbf{x}^i)^\top \textbf{C}^{-1} \textbf{y} = \textbf{C}^{i \top} \textbf{C}^{-1} \textbf{y} = \textbf{e}_i^\top \textbf{y} = y_i$.
This does not stand when $\textbf{C}^{-1}$ is replaced by $(\textbf{C}+\tau^2 \textbf{I})^{-1}$.

Recall that the term $\textbf{C}^{-1}\textbf{y}$ in the kriging mean of Equation~(\ref{E:kriging_mean}) is denoted by {\boldmath$\beta$}. 
When nugget regularization is used, {\boldmath$\beta$} is shown as {\boldmath$\beta$}$^{Nug}$ and, thanks to the eigenvalue decomposition of $(\textbf{C}+\tau^2 \textbf{I})^{-1}$, it is written
\begin{eqnarray} \label{beta_nug}
\bm{\beta}^{Nug} ~=~ \sum \limits_{i=1}^r \frac{\left(\textbf{V}^i \right)^\top \textbf{y}}{\lambda_i+\tau^2}\textbf{V}^i + \sum \limits_{i=r+1}^n \frac{\left(\textbf{W}^i \right)^\top \textbf{y}}{\tau^2}\textbf{W}^i.
\end{eqnarray}
The main difference between {\boldmath$\beta$}$^{PI}$ (Equation~(\ref{E:beta_PI})) and {\boldmath$\beta$}$^{Nug}$ lies in the second part of {\boldmath$\beta$}$^{Nug}$: the part that spans the null space of the covariance matrix. In the following, we show that this term cancels out when multiplied by $\textbf{c}(\textbf{x})^\top$, a product that intervenes in kriging.

\begin{property}[Orthogonality Property of \textbf{c} and Null(\textbf{C})] 
\mbox{~}\\
For all {\normalfont $\textbf{x} \in D$}, the covariance vector {\normalfont $\textbf{c}(\textbf{x})$} is perpendicular to the null space of the covariance matrix {\normalfont $\textbf{C}$}. 
\end{property}
\textit{Proof}: The kernel $K(.,.)$ is a covariance function \cite{aronszajn1950}, hence the matrix
\begin{eqnarray}
\textbf{C}_x=
\begin{bmatrix}
K(\textbf{x},\textbf{x}) & \textbf{c}(\textbf{x})^\top \\
\textbf{c}(\textbf{x}) & \textbf{C}
\end{bmatrix}
\end{eqnarray}
is positive semidefinite.

Let $\textbf{w}$ be a vector in the null space of $\textbf{C}$. According to the definition of positive semidefinite matrices, we have
\begin{eqnarray}
\begin{pmatrix}
1 \\
\textbf{w}
\end{pmatrix}^\top \textbf{C}_x
\begin{pmatrix}
1 \\
\textbf{w}
\end{pmatrix}=K(\textbf{x},\textbf{x}) + 2 \sum \limits_{i=1}^n K(\textbf{x}, x_i)w_i + 0 \geq 0.
\end{eqnarray}
The above equation is valid for any vector $\gamma \textbf{w}$ as well, in which  $\gamma$ is a real number. This happens only if $\sum \limits_{i=1}^n K(\textbf{x}, x_i)w_i$ is zero, that is to say,  $\textbf{c}(\textbf{x})^\top$ is perpendicular to the null space of  $\textbf{C}$. $\square$ \\

As a result of the Orthogonality Property of $\textbf{c}$ and Null(\textbf{C}), the second term of 
{\boldmath$\beta$}$^{Nug}$ in Equation~(\ref{beta_nug}) disappears in the kriging mean which becomes
\begin{eqnarray} \label{E:kriging_mean_nug}
m^{Nug}(\textbf{x})=\textbf{c}(\textbf{x})^\top \sum \limits_{i=1}^r \frac{\left(\textbf{V}^i \right)^\top \textbf{y}}{\lambda_i+\tau^2}\textbf{V}^i.
\end{eqnarray}
The Orthogonality Property applies similarly to the kriging covariance (Equation~(\ref{E:kriging_variance})), which yields 
\begin{equation} \label{E:kriging_variance_nug}
\begin{split}
c^{Nug}(\textbf{x},\textbf{x}^\prime) &= K(\textbf{x},\textbf{x}^\prime) - \textbf{c}(\textbf{x})^\top \sum \limits_{i=1}^r \frac{\left(\textbf{V}^i\right)^\top \textbf{c}(\textbf{x}^\prime)}{\lambda_i + \tau^2} \textbf{V}^i \\
& =K(\textbf{x},\textbf{x}^\prime) - \sum \limits_{i=1}^r \frac{\left(\left(\textbf{V}^i\right)^\top \textbf{c}(\textbf{x})\right) \left(\left(\textbf{V}^i\right)^\top \textbf{c}(\textbf{x}^\prime)\right) }{\lambda_i + \tau^2}.
\end{split}
\end{equation}
Comparing Equations (\ref{E:kriging_mean_PI}) and (\ref{E:kriging_mean_nug}) indicates that the behavior of $m^{PI}$ and $m^{Nug}$ will be similar to each other if $\tau^2$ is small. 
The same holds for kriging covariances (hence variances) $c^{PI}$ and $c^{Nug}$ in Equations (\ref{E:kriging_variance_PI}) and (\ref{E:kriging_variance_nug}). 
\begin{property}[Equivalence of PI and nugget regularizations] 
\label{prop-PInug_equiv}
\mbox{~}\\
The mean and covariance of conditional GPs regularized by nugget tend toward the ones of GPs regularized by pseudoinverse as the nugget value $\tau^2$ tends to 0.
\end{property}
In addition, Equations (\ref{E:kriging_variance_PI}) and (\ref{E:kriging_variance_nug}) show that $c^{Nug}$ is always greater than  $c^{PI}$. These results will be illustrated later in the Discussion Section. 

\subsection{Nugget and maximum likelihood} \label{nugget_and_ML}
It is common to estimate the nugget parameter by maximum likekihood (ML, cf. Appendix~\ref{sec-appnugget}, Equation~(\ref{E:nugget_likelihood})). 
As will be detailed below, the amplitude of the nugget estimated by ML is increasing with the spread of observations at redundant points. 
It matches the interpretation of nugget as the amount of noise put on data: an increasing discrepancy between responses at a given point is associated to more observations noise.

\begin{figure}[htpb]
\begin{center}
\includegraphics[width=0.5\textwidth]{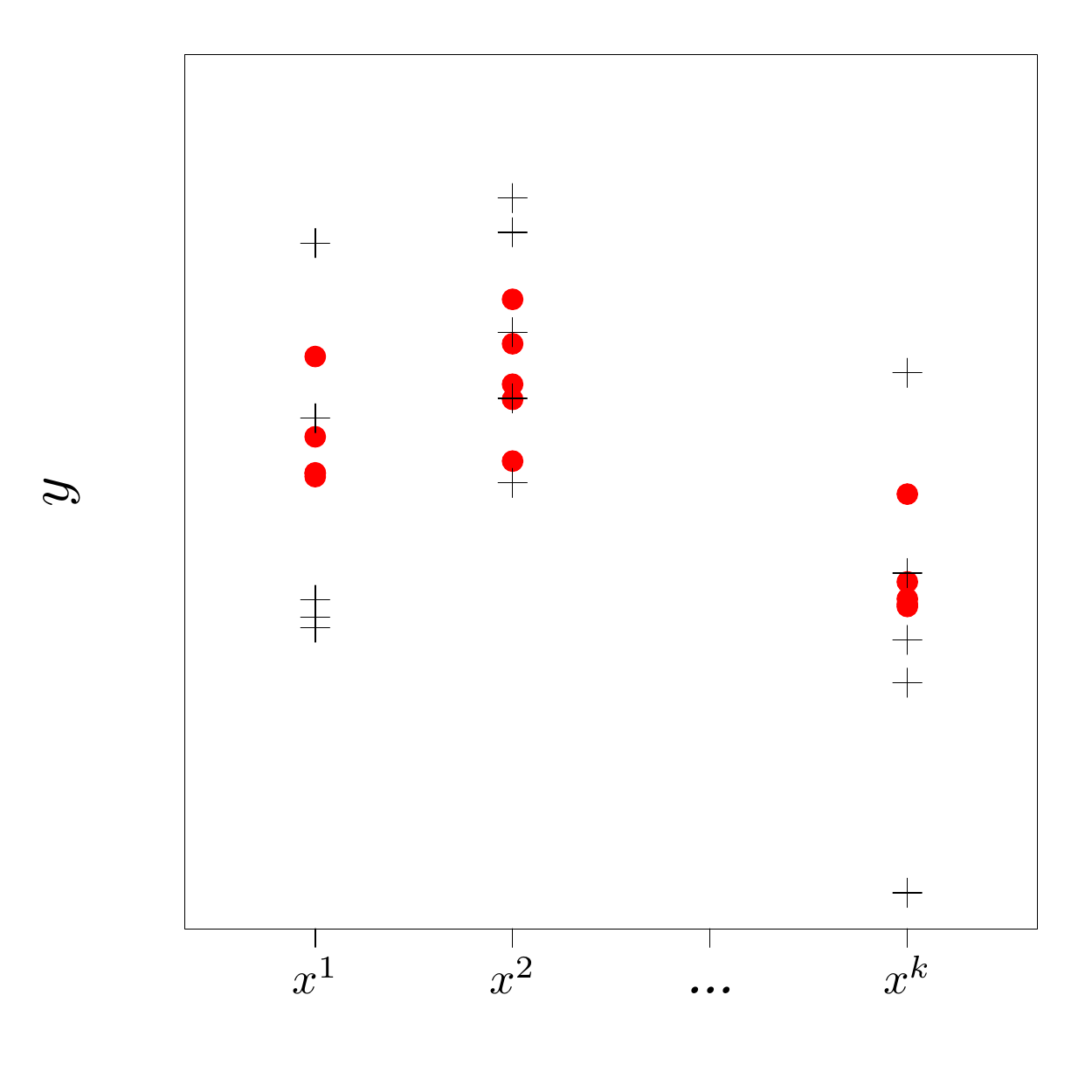}
   \caption{
The response values $\textbf{y}$ and $\textbf{y}^+$ are denoted by bullets and crosses, respectively. 
At each location, the mean of $\textbf{y}$ and $\textbf{y}^+$ are identical, $\overline{y_i}=\overline{{y^+}_i}$, 
but the spread of observations in $\textbf{y}^+$ is never less than that of $\textbf{y}$ at redundant points.
\label{F:redundant_points}
}
\end{center}
\end{figure}

In Figure~\ref{F:redundant_points} two vectors of response values are shown, $\textbf{y}$ (bullets) and $\textbf{y}^+$ (crosses), located at $k$ different $\textbf{x}$ sites.
The spread of response values $\textbf{y}^+$ is larger than that of $\textbf{y}$ at some redundant points. 
Let $s_i^2$ and ${s_i^+}^2$, \linebreak $i=1, ..., k$, denote the variances of $\textbf{y}$ and $\textbf{y}^+$ at the redundant points,
\begin{equation}\label{var_redundant}
s_i^2=\frac{\sum \limits_{j=N_1+...+N_{i-1}+1}^{N_1+...+N_i}\left( y_j - \overline{y_i}\right)^2}{N_i - 1}, 
\end{equation}
and the same stands with $\textbf{y}^+$ and its variance ${s_i^+}^2$.

The nugget that maximizes the likelihood, the other GP parameters being fixed (the length-scales $\theta_i$  and the process variance $\sigma^2$),
is increasing when the variance of the outputs increases.
\newtheorem{theorem}{Theorem}
\begin{theorem}
 \mbox{~}\\
Suppose that there are observations located at {\normalfont $k$} different sites. 
If we are given two vectors of response values {\normalfont $\textbf{y}$} and {\normalfont $\textbf{y}^+$} such that
\begin{enumerate}
\item {\normalfont ${s_i^+}^2\geq {s_i^2}$} for all {\normalfont $i=1, \ldots, k $} and 
\item {\normalfont $\overline{y_i}=\overline{{y^+}_i}~$} for all {\normalfont $i=1, ..., k$},  
\end{enumerate}
then the nugget amplitudes  {\normalfont $\hat{\tau}^2$} and {\normalfont $\widehat{\tau^+}^2$} 
that maximize the likelihood with other GP parameters being fixed are 
such that {\normalfont $\widehat{\tau^+}^2 \geq \hat{\tau}^2$}.
\label{theorem-MLnugget}
\end{theorem} 
A proof for Theorem \ref{theorem-MLnugget} is given in Appendix~\ref{sec-appnugget}.

\section{Discussion: choice and tuning of the classical regularization methods} \label{discussion}
This section carries out a practical comparison of PI and nugget regularization methods, which are readily available in most GP softwares \cite{siefert, roustant2012}.
We start with a discussion of how data and model match, which further allows to decide whether nugget or PI should be used. Finally, we provide guidelines to tune the regularization parameters.

Before starting with our discussion, note that nugget regularization is the method of choice 
when the observed data is known to be corrupted by an additive noise that is homogeneous in $D$ 
since, in this case, the nugget amplitude $\tau^2$ has the physical meaning of noise variance \cite{roustant2012}. 
The loss of the interpolating property at data points associated to nugget regularization is 
then an intended filtering effect. 
The rest of our discussion will assume no knowledge of the eventual noise model affecting observations.

\subsection{Model-data discrepancy}
Model-data discrepancy can be measured as the distance between the observations $\textbf y$ and the GP model regularized by pseudoinverse.
\begin{definition}[Model-data discrepancy]\label{def-discrepancy}
Let {\normalfont $\mathbf X$} be a set of design points with associated observations {\normalfont $\mathbf{y}$}. Let {\normalfont $\mathbf V$} and {\normalfont $\mathbf{W}$} be the normalized eigenvectors spanning the image space and the null space of the covariance matrix {\normalfont $\mathbf{C}$}, respectively. The model-data discrepancy is defined as 
\begin{equation}
{\normalfont discr ~=~\frac{\Vert \mathbf y - m^{PI}(\mathbf X) \Vert^2}{\Vert \mathbf y \Vert^2} 
~=~ \frac{\Vert \textbf{WW}^\top \mathbf y \Vert^2}{\Vert \mathbf y \Vert^2} }
\label{eq-discr}
\end{equation}
where {\normalfont $m^{PI}(\ldotp)$} is the pseudoinverse regularized GP model of Equation~(\ref{eq-mpi}).
\end{definition}
The last equality in the definition of $discr$ comes from Equations~(\ref{eq-PIVVT}) and (\ref{eq-Pim}).
The discrepancy is a normalized scalar, $0 \le discr \le 1$, where $discr=0$ indicates that the model and the data are perfectly compatible, and vice versa when $discr=1$.

The definition of redundant points does not depend on the observations $\textbf y$ and the model-data discrepancy is a scalar globalizing the contributions of all observations.  An intermediate object between redundant points and discrepancy is the gradient of the squared model-data error with respect to the observations,
\begin{equation}
\nabla_{\textbf y} \Vert \textbf y - m^{PI}(\textbf X) \Vert^2 ~=~ \textbf{WW}^\top \textbf y ~.
\label{eq-graderror}
\end{equation}
It appears that the gradient of the error, $\Vert \textbf y - m^{PI}(\textbf X) \Vert^2$, is equal to the model-data distance, $\textbf{W}\textbf{W}^\top\textbf{y}$. This property comes from the quadratic form of the error. 
The magnitude of the components of the vector $\textbf{WW}^\top \textbf y$ measure the sensitivity of the error to a particular observation. At repeated points, a gradient-based approach where the $y$'s are optimized would advocate to make the observations closer to their mean proportionally to their distance to the mean.  
In other words, $- \textbf{WW}^\top \textbf y$ is a direction of reduction of the model-data distance in the space of observations.
Because the distance considered is quadratic, this direction is colinear to the error, ~($\textbf y - m^{PI}(\textbf X)$).
The indices of the non-zero components of $\textbf{WW}^\top \textbf y$ also designate the redundant points.

\subsection{Two detailed examples}
\begin{figure}[htpb] 
\centering
\includegraphics[width=0.5\textwidth]{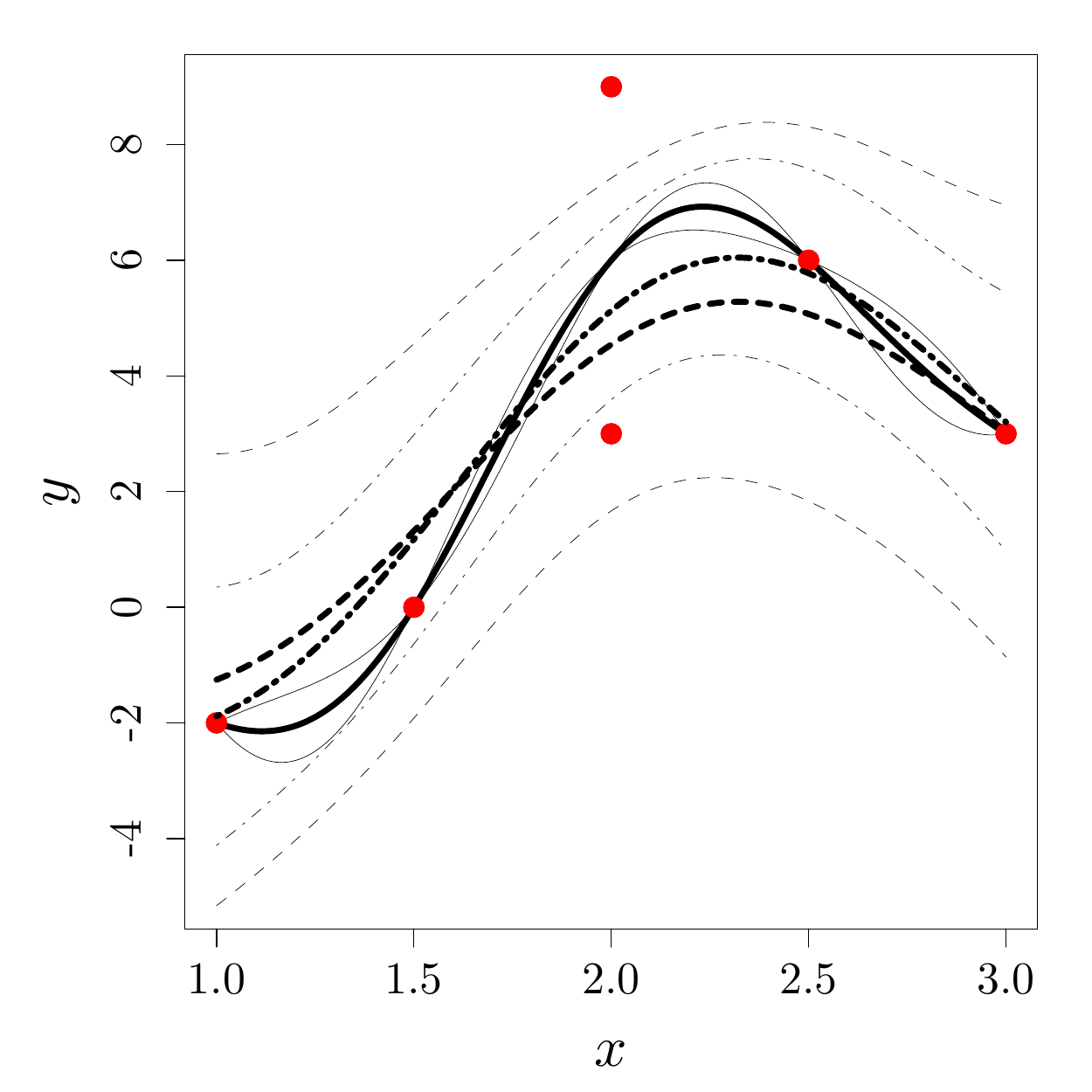}
   \caption{
Comparison of kriging regularized by PI (solid lines), nugget estimated by ML (dashed lines) and nugget estimated by cross-validation (dash-dotted lines). 
$\textbf X = [1;1.5;2;2.00001;2.5;3]$ and $\textbf y = (-2,0,3,9,6,3)^\top$.
The estimated nugget values are $\hat{\tau}^2=7.06$ and  $\hat{\tau}^2_{CV}=1.75$.  
\label{nugget_estim}
}
\end{figure}

A common practice when the nugget value, $\tau^2$, is not known beforehand is to estimate it by ML or cross-validation. 
In Appendix \ref{sec-appnugget}, we show that the ML estimated nugget value, $\hat{\tau}^2$, 
is increasing with the spread of responses at redundant points. 
This is one situation (among others, e.g., the additive example hereafter) where the data and the model mismatch, and $\hat{\tau}^2$ is large.
Figure~\ref{nugget_estim} is an example where $\hat{\tau}^2$ is equal to 7.06. The kriging mean and the 68\% confidence interval, in this case of nugget regularization tuned by ML, are drawn with dashed lines. 
Some authors such as in \cite{wagner2010, bachoc2013} recommend using cross-validation instead of ML for learning the kriging parameters. 
In the example of Figure~\ref{nugget_estim}, the estimated nugget value by leave-one-out cross-validation, denoted by $\hat{\tau}^2_{CV}$, is 1.75.
The dash-dotted lines represent the kriging model regularized by nugget that is estimated by cross-validation. 
The model-data discrepancy is $discr=0.36$ and $\textbf{WW}^\top \textbf y = (0,0,-3,3,0,0)^\top $ which shows that points 3 and 4 are 
redundant and their outputs should be made closer to reduce the model-data error. Whether or not in practice the outputs can be controlled is out of the scope of our discussion. 

\begin{figure}[htpb]
\centering
\includegraphics[width=0.49\textwidth]{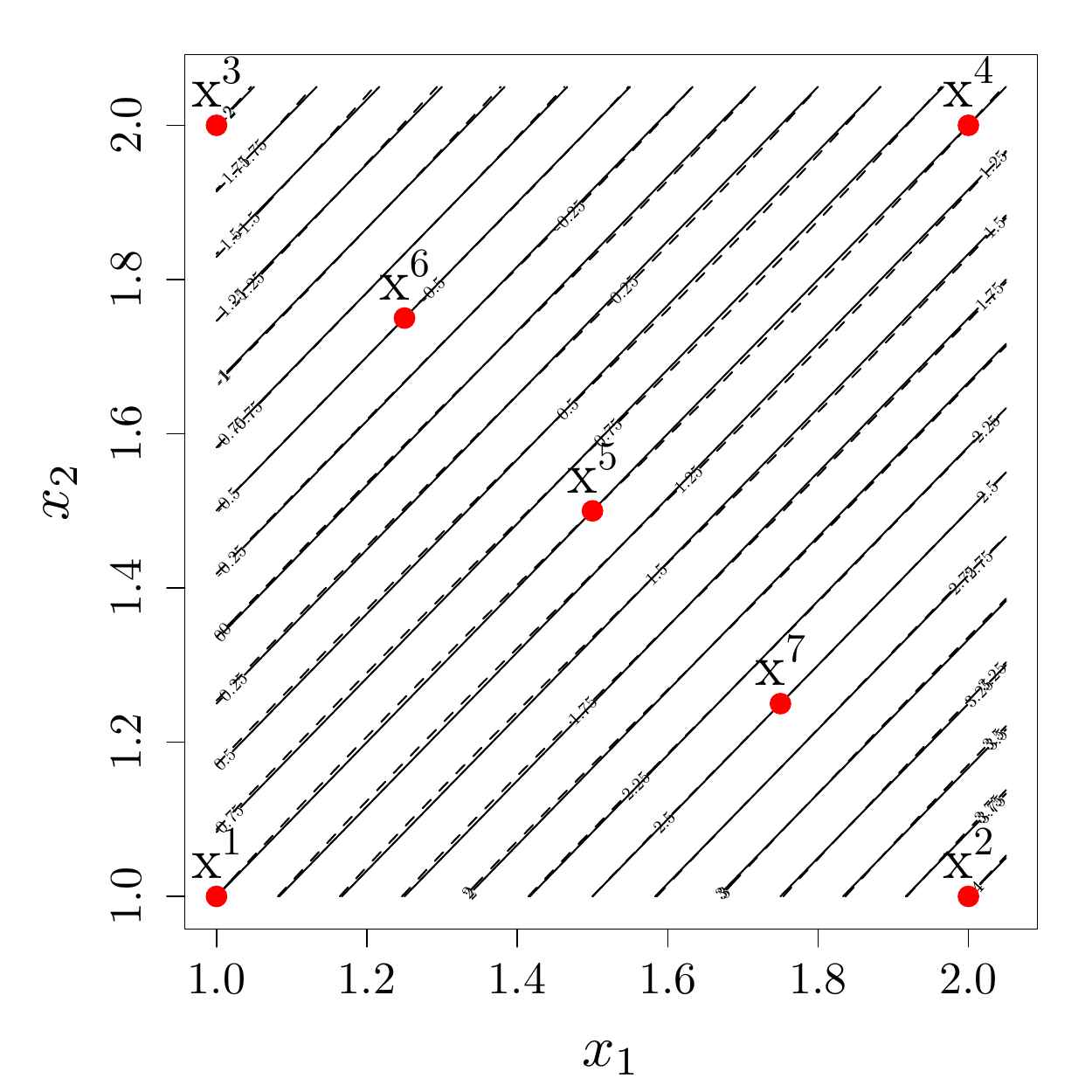}
\includegraphics[width=0.49\textwidth]{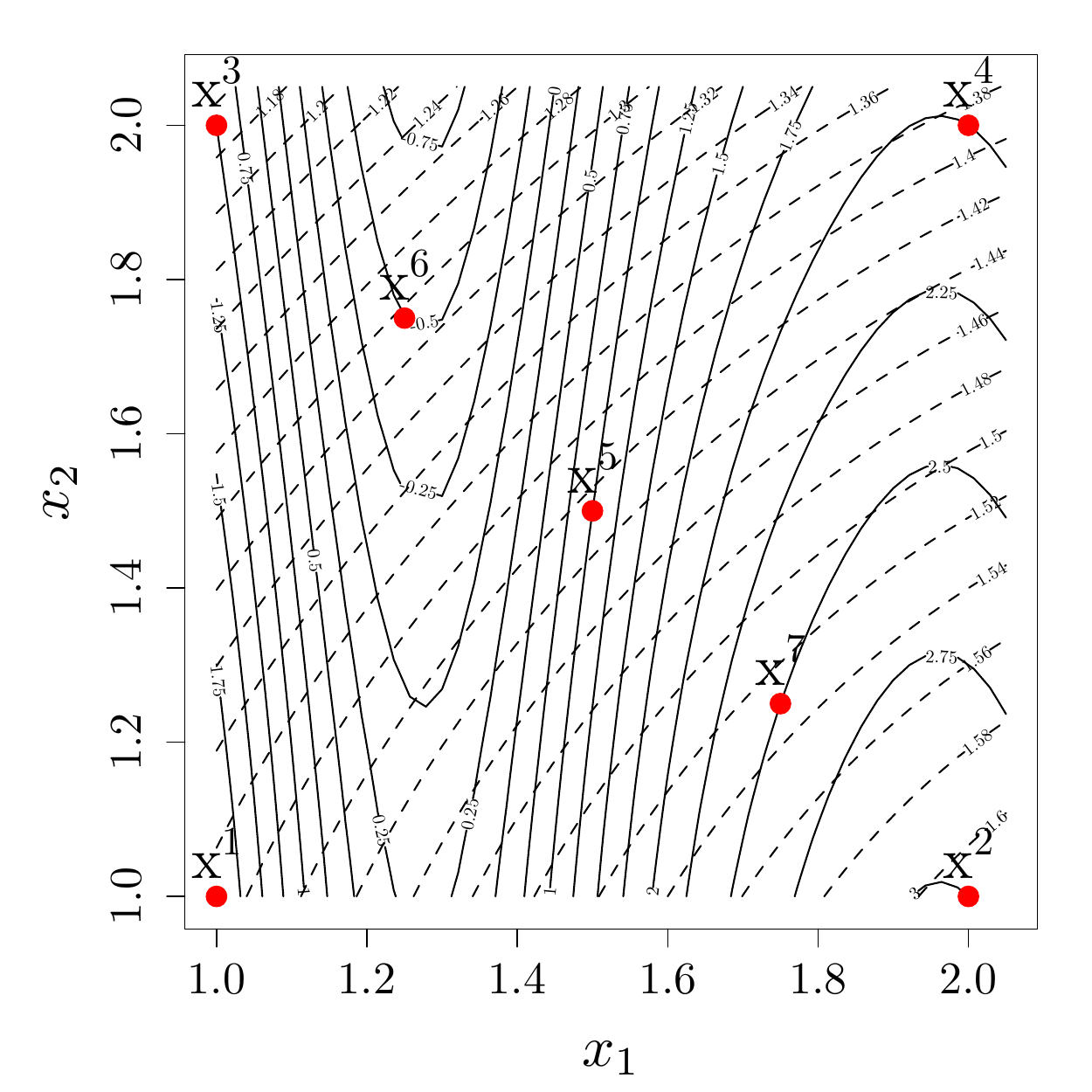}
\caption{Contour plots of kriging mean regularized by pseudoinverse (solid line) vs. nugget (dashed line) for an additive GP. 
The bullets are data points. Left: the response values are additive, $\textbf{y}=(1, 4, -2, 1, 1, -0.5, 2.5)^\top$ and 
$\hat{\tau}^2 = 10^{-12}$. 
Right: the third observation is replaced by 2, creating non-additive observations and $\hat{\tau}^2 \approx 1.91$; $m^{Nug}(\textbf{x})$ 
is no longer interpolating, $m^{PI}(\textbf{x})$ still interpolates $\textbf{x}^5$ to $\textbf{x}^7$.}
\label{additive_kernel}
\end{figure}

We now give a two-dimensional example of a kriging model with additive kernel defined over 
$\textbf{X}=\left[(1, 1), (2, 1), (1, 2), (2, 2), (1.5, 1.5), (1.25, 1.75), (1.75, 1.25) \right]$, cf. Figure \ref{additive_kernel}. 
As explained in Section~\ref{degeneracy_cov.mat}, the first four points of the DoE make the additive covariance matrix non-invertible 
even though the points are not near each other in Euclidean distance.
Suppose that the design points have the response values $\textbf{y}=(1, 4, -2, 1, 1, -0.5, 2.5)^\top$ which correspond to the additive true function $f(\textbf{x})=x^2_1-x^2_2 + 1$. The covariance matrix is the sum of two parts
\begin{equation*}
\textbf{C}_{add}=\sigma^2_1 K_1 + \sigma^2_2 K_2~,
\end{equation*}
where $\sigma^2_i$ are the process variances and $\sigma^2_iK_i$ the kernel in dimension $i=1, 2$. 

To estimate the parameters of $\textbf{C}_{add}$, the negative of the likelihood is numerically 
minimized (see Equation~(\ref{E:nugget_likelihood})) which yields a nugget value 
$\hat{\tau}^2 \approx 10^{-12}$ (the lower bound on nugget used).
A small nugget value is obtained because the associated output value follows an additive function compatible with the kernel:
there is no discrepancy between the model and the data.
Because of the small nugget value, the models regularized by PI and nugget are very close to each other (the left picture in Figure \ref{additive_kernel}).

We now introduce model-data discrepancy by changing the third response from -2 to 2: 
additive kriging models can no longer interpolate these outputs.
The nugget value estimated by ML is equal to $1.91$, so $m^{Nug}(\textbf{x})$ does not interpolate any of the data points 
($\textbf{x}^1$ to $\textbf{x}^7$). 
Regarding $m^{PI}(\textbf{x})$, the projection onto $Im(\textbf C)$ make the GP predictions different from the observations at $\textbf{x}^1$,
$\textbf{x}^2$, $\textbf{x}^3$ and $\textbf{x}^4$. For example, $m^{PI}(\textbf{x}^4)=2$ when $\textbf{y}_4=1$.
The projection applied to points $\textbf{x}^5$ to $\textbf{x}^7$ (where no linear dependency between the associated columns of $\textbf{C}$ exists) show that $m^{PI}(\textbf{x})$ is interpolating 
there, which is observed on the right picture of Figure~\ref{additive_kernel}.

The above observations point out that large estimated values of nugget (whether by ML or cross-validation) indicate model-data discrepancy, in agreement with the calculated discrepancies of Eq.~(\ref{eq-discr}): in the last additive kernel example when all the outputs were additive, $discr=0$ and $\textbf{WW}^\top \textbf y = (0,0,0,0,0,0,0)^\top$ (no redundant point); when the value of the third output was increased to 2, $discr=0.37$ and $\textbf{WW}^\top \textbf y = (-1,1,1,-1,0,0,0)^\top$ showing that points 1 to 4 are redundant and that, to reduce model error, points 1 and 4 should increase their outputs 
while points 2 and 3 should decrease theirs. 

For the sole purpose of quantifying model-data discrepancy, it is more efficient to calculate $discr$ 
than using the estimated nugget. Formula~(\ref{eq-discr}) involves
one pseudo-inverse calculation and two matrix products when nugget estimation implies 
a nonlinear likelihood maximization with repeated embedded $\textbf{C}$ eigenvalues analyses.

\subsection{PI or nugget ?}
\begin{figure}[htpb]
\centering
\includegraphics[width=0.49\textwidth]{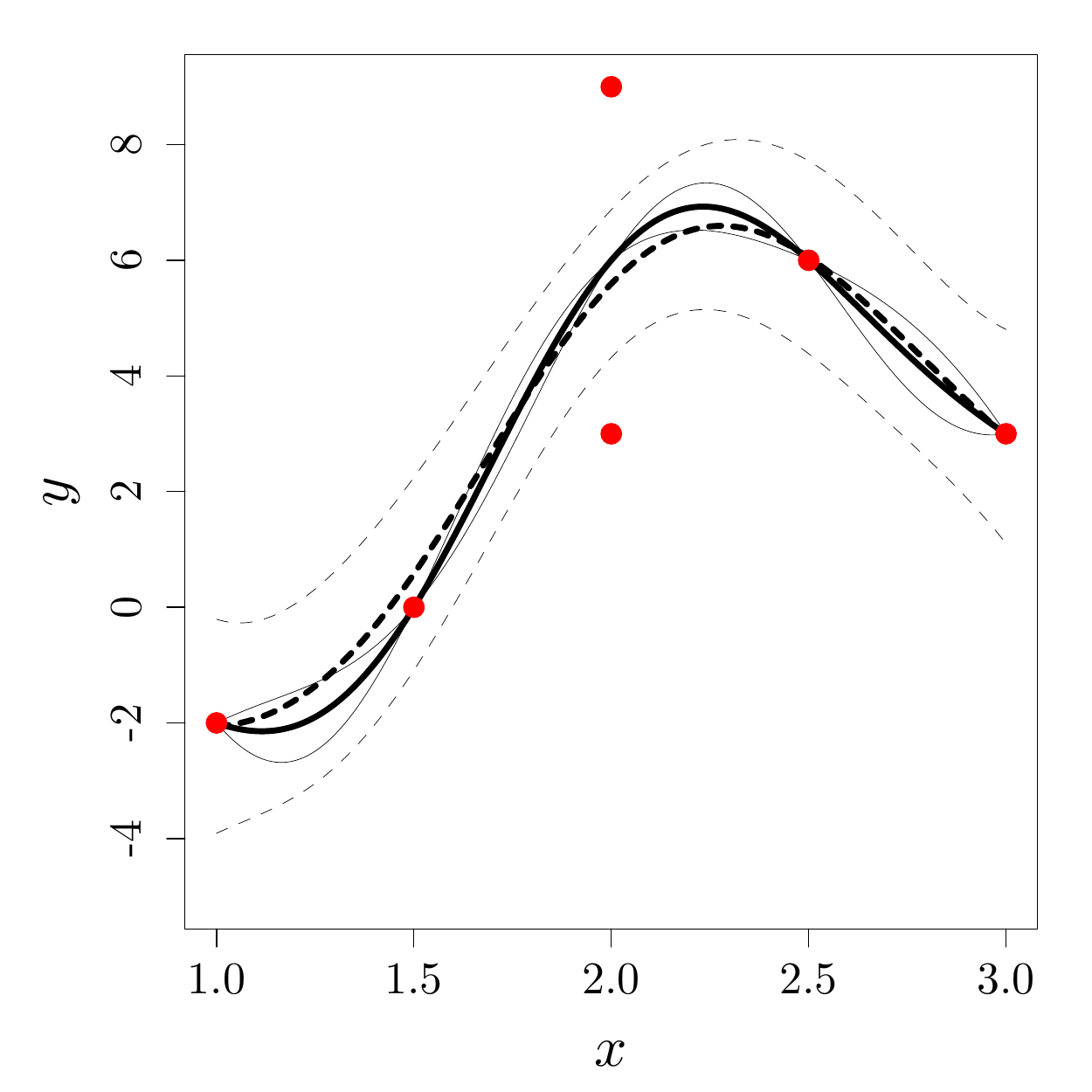}
\includegraphics[width=0.49\textwidth]{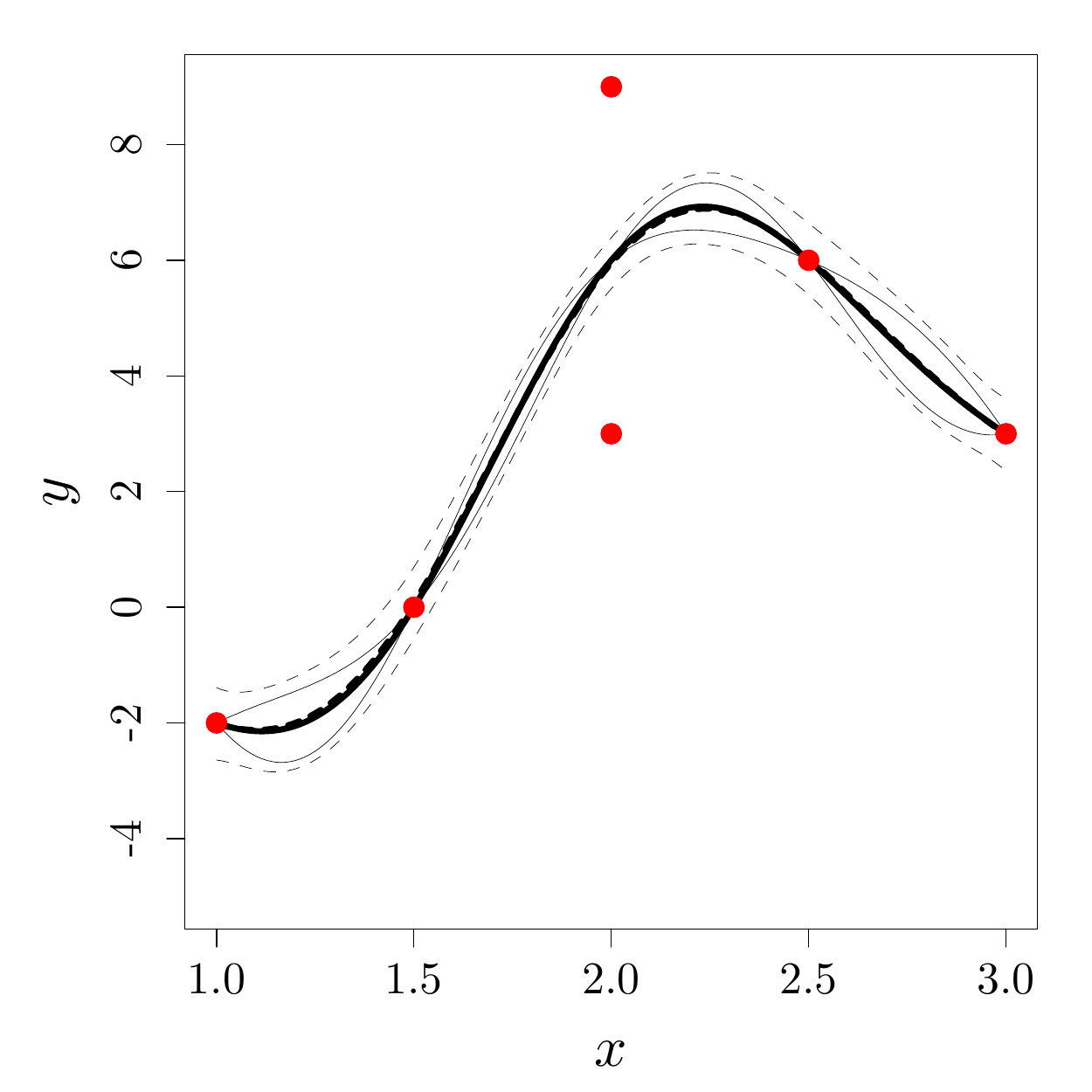}
\vspace{-0.1cm}
   \caption{
One dimensional kriging regularized by PI (solid lines) and nugget (dashed lines). 
The nugget amplitude is 1 on the left and 0.1 on the right. The cut-off eigenvalue for the pseudoinverse 
(below which eigenvalues are rounded off to 0) is $\eta=10^{-3}$. $m^{Nug}(x)$ is not
interpolating which is best seen at the second point on the left. On the right, the PI and nugget models are closer to each other. 
Same $\textbf X$ and $\textbf y$ as in Figure~\ref{nugget_estim}.
\label{fixed_nugget}
}
\end{figure}

On the one hand, models regularized by PI have predictions, $m^{PI}(\ldotp)$, that interpolate uniquely defined points and 
go through the average output at redundant points (Property \ref{PI-averaging}). 
The associated kriging variances, $v^{PI}(\ldotp)$, are null at redundant points (Property \ref{PI-nulvar}).
On the other hand, models regularized by nugget have predictions which are neither interpolating nor averaging (Property \ref{prop-lossinterp}) 
while their variances are non-zero at data points. 
Note that kriging variance tends to $\sigma^2$ as the nugget value increases (see Equation~(\ref{E:kriging_variance_nug})). 
These facts can be observed in Figure~\ref{fixed_nugget}. Additionally, this 
Figure illustrates that nugget regularization tends to PI regularization as the nugget value decreases (Property~\ref{prop-PInug_equiv}).
If there is a good agreement between the data and the GP model, the PI regularization or equivalently, a small nugget, should be used.
This can also be understood through the Definition of model-data discrepancy and Property~\ref{prop-PIproj}: when $discr=0$, the observations 
are perpendicular to $Null(\textbf{C}) $ and, equivalently, $m^{PI}(\textbf{X})=\textbf{y}$ since $m^{PI}(\ldotp)$ performs a 
projection onto $Im(\textbf{C})$.
Vice versa, if the model-data discrepancy measure is significant, 
choosing PI or nugget regularization will have a strong impact on the model: 
either the prediction averaging property is regarded as most important and PI should be used, or a non-zero variance at redundant points 
is favored and nugget should be selected; If the discrepancy is concentrated on few redundant points, nugget regularized models will 
distribute the uncertainty (additional model variance) throughout the $x$ domain while PI regularized models will ignore it. 
From the above arguments, it is seen that the decision for using PI or nugget regularizations is problem dependent.

\subsection{Tuning regularization parameters}
Nugget values may be estimated by maximum likelihood (\cite{roustant2012}) or cross-validation (\cite{cressie1993,bachoc2013}), but as argued in the previous Section, it may also be preferred to fix them to small numbers.
How small can a nugget value be? Adding nugget to the main diagonal of a covariance matrix increments all the eigenvalues by the nugget amplitude. The condition number of the covariance matrix with nugget is $\kappa(\textbf{C} + \tau^2 \textbf{I})=\frac{\lambda_{max}+\tau^2}{\lambda_{min}+\tau^2}$. Accordingly, a ``small" nugget is the smallest value of $\tau^2$ such that $\kappa(\textbf{C} + \tau^2 \textbf{I})$ 
is less than a reasonable condition number after regularization, $\kappa_{max}$ (say, $\kappa_{max} = 10^8$).
With such targeted condition number, the smallest nugget would be $ \tau^2 ~=~ \frac{\lambda_{max} - \kappa_{max} \lambda_{min}}{\kappa_{max}-1} $ ~if~ $ \lambda_{max} - \kappa_{max} \lambda_{min} \ge 0$,  $ \tau^2 ~=~ 0$ ~ otherwise.

Computing a pseudoinverse also involves a parameter, the positive threshold $\eta$ below which an eigenvalue is considered as null. 
The eigenvectors associated to eigenvalues smaller than $\eta$ are numerically regarded as null space basis vectors (even though they 
may not, strictly speaking, be part of the null space).
A suitable threshold should filter out eigenvectors associated to points that are almost redundant. The heuristic we propose is to tune $\eta$ so that 
$\lambda_1/\eta$, which is an upper bound of the PI condition number\footnote{
By PI condition number we mean $\kappa_{PI}(\textbf{C}) = \| \textbf{C} \| \| \textbf{C}^\dagger \| = \lambda_1/\lambda_r \le \lambda_1/\eta  $
}, is equal to $\kappa_{max}$, i.e., $\eta =  \lambda_1 / \kappa_{max} $.

\begin{figure}[htpb]
\centering
\includegraphics[width=0.5\textwidth]{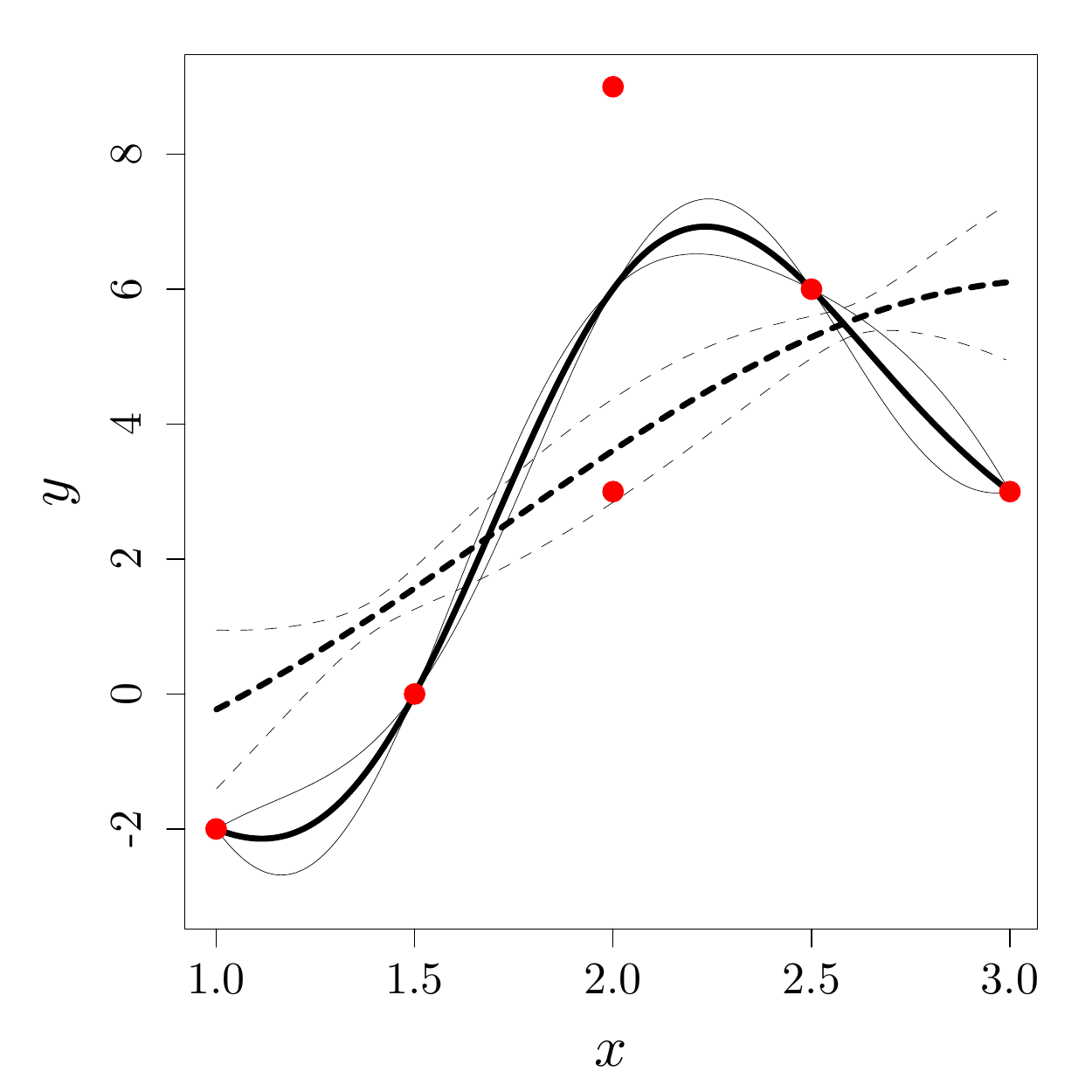}
   \caption{
Effect of the tolerance $\eta$ on the kriging model regularized by PI.
Dashed line, $\eta=1$; continuous line, $\eta=10^{-3}$.
Except for $\eta$, the setting is the same as that of Figure \ref{fixed_nugget}.
When the tolerance is large ($\eta=1$), the 5th eigenvector is deleted from the effective image space of $\textbf{C}$ in addition to the 6th eigenvector, and the PI regularized model is no longer interpolating.
Same $\textbf X$ and $\textbf y$ as in Figure~\ref{nugget_estim}.
\label{PI_tau}
}
\end{figure}

In the example shown in Figure \ref{PI_tau}, the covariance matrix is not numerically invertible because the points 3 and 4 are near $x=2$. 
The covariance matrix has six eigenvalues, $\lambda_1 = 34.89 \geq ... \geq \lambda_5 = 0.86 \geq \lambda_6 = 8.42\times 10^{-11} \approx 0$ and the eigenvector related to the smallest 
eigenvalue is $\textbf{W}^1 = (\textbf{e}^4 - \textbf{e}^3)/\sqrt{2}$. 
In Figure~\ref{PI_tau}, we have selected $\eta = 10^{-3}$, hence $\kappa_{PI}(\textbf{C})=40.56$. 
Any value of $\eta$ in the interval $\lambda_6 < \eta < \lambda_5$ would have yielded the same result.
But if the selected tolerance were e.g., $\eta=1$, which is larger than $\lambda_5$, the obtained PI kriging model no longer interpolates
data points. 

\section{Interpolating Gaussian distributions} 
\label{sec-interpol}

\subsection{Interpolation and repeated points}
In our context of deterministic experiments, we are interested in interpolating data.
The notion of interpolation should be clarified in the case of repeated points with 
different outputs (e.g., Figure~\ref{F:redundant_points}) as a function cannot interpolate them. 
Interpolation should then be generalized to redundant points, which comprise repeated points, 
but to keep the explanations simple we do not do it in the current paper.
Here, we seek GPs that have the following interpolation properties.
\begin{definition}[Interpolating Gaussian Process]
A GP interpolates a given set of data points when
\begin{itemize}
\item its trajectories pass through uniquely defined data points (therefore the GP has a null variance there),
\item and at repeated points the GP's mean and variance are the empirical average and variance of the outputs, respectively. 
\end{itemize}
\end{definition}
The following GP model has the above interpolation properties for deterministic outputs, even in the presence of repeated points. 
In this sense, it can be seen as a new regularization technique, although its potential use goes beyond regularization.

\subsection{A GP model with interpolation properties}
We now introduce a new GP model, called \textit{distribution-wise} GP, with the desirable interpolation properties in the presence of repeated points. 
Accordingly, it can be regarded as a regularization method. 
Moreover, it is computationally more efficient than the \textit{point-wise} GP models.

Following the same notations as in Section~\ref{PI_averaging}, the model is built from observations at $k$ \emph{different} $\textbf{x}$ sites. 
The basic assumption that makes distribution-wise GP different from point-wise GP is that observations 
are seen as realizations of a known joint Gaussian probability distribution.
In distribution-wise GP, it is assumed that a distribution is observed at each location, 
as opposed to usual conditional GPs where unique values of the process are observed. 
This is the reason for the name ``distribution-wise GP''.
Observing distributions brings in a framework which, as will now be seen, 
is compatible with repeated points.

Let $Z(\textbf{x}^i) \sim \mathcal{N}(\mu_{Z_i}, \sigma^2_{Z_i})$ denote the random outputs at locations $\textbf{x}^i~,~i=1,\dots, k~$ and their known probability distributions. 
Together, the $k$ sets of observations make the random vector $\textbf{Z}=(Z(\textbf{x}^1), \ldots, Z(\textbf{x}^k)) \sim \mathcal{N}(\boldsymbol\mu_{Z}, \boldsymbol\Gamma_{Z})$ in which the diagonal of the matrix $\boldsymbol\Gamma_{Z}$ is made of the $\sigma^2_{Z_i}$'s.

The distribution-wise GP is derived in two steps through conditioning: first it is assumed that the vector $\textbf{Z}$ is given, and the usual conditional GP (kriging) formula can be applied; then the randomness of $\textbf{Z}$ is accounted for. 
The conditional mean and variance of the distribution-wise GP, $m^{Dist}$ and $v^{Dist}$, come from the laws of total expectation and variance applied to $\textbf{Z}$ and the GP outcomes $\omega \in \Omega$:
\begin{align}
\nonumber m^{Dist}(\textbf{x})&~=\mathbb{E}_Z \left( \mathbb{E}_\Omega (Y(\textbf{x})|Y(\textbf{x}^i)= Z(\textbf{x}^i) ~,~ 1 \leq i \leq k  \right)~=~ \\
& \mathbb{E}_Z \left( \textbf{c}_Z(\textbf{x})^\top \textbf{C}_Z^{-1} \textbf{Z} \right )
~=~ \textbf{c}_Z(\textbf{x})^\top \textbf{C}_Z^{-1} \boldsymbol\mu_{Z} 
\end{align}
where the $Z$ subscript is used to distinguish between the point-wise and the distribution-wise covariances. $\textbf{C}$ is $n \times n$ and not necessarily invertible while $\textbf{C}_Z$ is $k \times k$ and invertible because the $k$ $\textbf{x}^i$'s are different\footnote{Remember that only repeated points are considered here. For the more general redundant points, invertibility of $\textbf{C}_Z$ will be needed but the way to achieve it is out of the scope of this paper.}.
The variance is calculated in a similar way
\begin{align}
\nonumber v^{Dist}(\textbf{x})&=~\mathbb{E}_Z \left( \Var_\Omega (Y(\textbf{x}) | Y(\textbf{x}^i)= Z(\textbf{x}^i) ~,~ 1 \leq i \leq k  \right)+ \\ \nonumber & \Var_Z \left( \mathbb{E}_\Omega (Y(\textbf{x}) | Y(\textbf{x}^i)= Z(\textbf{x}^i) ~,~ 1 \leq i \leq k  \right)=~ \\
&\textbf{c}_Z(\textbf{x},\textbf{x}) - \textbf{c}_Z(\textbf{x})^\top \textbf{C}_Z^{-1}\textbf{c}_Z(\textbf{x})~+~ \textbf{c}_Z(\textbf{x})^\top \textbf{C}_Z^{-1} \underbrace{(\Var_Z\textbf{Z})}_{\boldsymbol\Gamma_{Z}} \textbf{C}_Z^{-1}\textbf{c}_Z(\textbf{x}) ~.
\end{align}

The distribution-wise GP model interpolates the mean and the variance of the distributions at the $k$ locations. At an arbitrary location $i$, the term $\textbf{c}_Z(\textbf{x})^\top \textbf{C}_Z^{-1}$ that appears in both $m^{Dist}$ and $v^{Dist}$ becomes $\textbf{e}^\top_i$ because $\textbf{c}_Z(\textbf{x}^i)$ is the $i$th column of $\textbf{C}_Z$ in this case. As a result 
\begin{align}
& m^{Dist}(\textbf{x}^i) = \textbf{c}_Z(\textbf{x}^i)^\top \textbf{C}_Z^{-1} \boldsymbol\mu_{Z}=\mu_{Z_i} \qquad \text{ and }\\
\nonumber  &v^{Dist}(\textbf{x}^i) = \textbf{c}_Z(\textbf{x}^i,\textbf{x}^i) - \textbf{c}_Z(\textbf{x}^i)^\top \textbf{C}_Z^{-1}\textbf{c}_Z(\textbf{x}^i)~+~ \\ 
&\qquad \qquad \qquad \qquad \textbf{c}_Z(\textbf{x}^i)^\top \textbf{C}_Z^{-1} \boldsymbol\Gamma_{Z} \textbf{C}_Z^{-1}\textbf{c}_Z(\textbf{x}^i) = \sigma^2_{Z_i}.
\end{align}

In practice, $\boldsymbol\mu_{Z}$ and $\boldsymbol\Gamma_{Z}$ can be approximated by the empirical mean and variance. 
Recall repeated points are grouped by sites, e.g., $y_1 , \ldots , y_{N_1}$ are the observations at $\textbf{x}^1$. The output empirical mean and variance at $\textbf{x}^i$ are $\overline{y_i}$ and $\overline{s^2_i}$ that we gather in the vector $\overline{\textbf{y}}$ and the $k \times k$ diagonal matrix $\hat{\textbf{$\Gamma$}}$ made of the $\overline{s^2_i}$'s. Then, the mean and the variance of the distribution-wise GP are expressed as 
\begin{align}
m^{Dist}(\textbf{x}) &~ \equiv ~\textbf{c}_Z(\textbf{x})^\top \textbf{C}_Z^{-1} \overline{\textbf{y}},\\
v^{Dist}(\textbf{x}) &~\equiv~\textbf{c}_Z(\textbf{x},\textbf{x}) - \textbf{c}_Z(\textbf{x})^\top \textbf{C}_Z^{-1}\textbf{c}_Z(\textbf{x})~+~ \textbf{c}_Z(\textbf{x})^\top \textbf{C}_Z^{-1} \hat{\textbf{$\Gamma$}}\textbf{C}_Z^{-1}\textbf{c}_Z(\textbf{x}). 
\end{align}
Because detecting repeated points has a computational cost of at most $\mathcal O(n^2)$, 
the computational complexity of the distribution-wise GP is driven by the covariance matrix inversion, 
as is the case for point-wise GP: distribution-wise GP applied to repeated points scales
in $\mathcal O(k^3)$ when traditional GP scale in $\mathcal O(n^3)$. I
The advantage of distribution-wise GP grows with the number of repeated points.
This result might not stand for 
the most general redundant points whose characterization may imply an $\mathcal O(n^3)$ eigendecomposition of $\textbf{C}$.

\begin{figure}[htpb]
\centering
\includegraphics[width=0.5\textwidth]{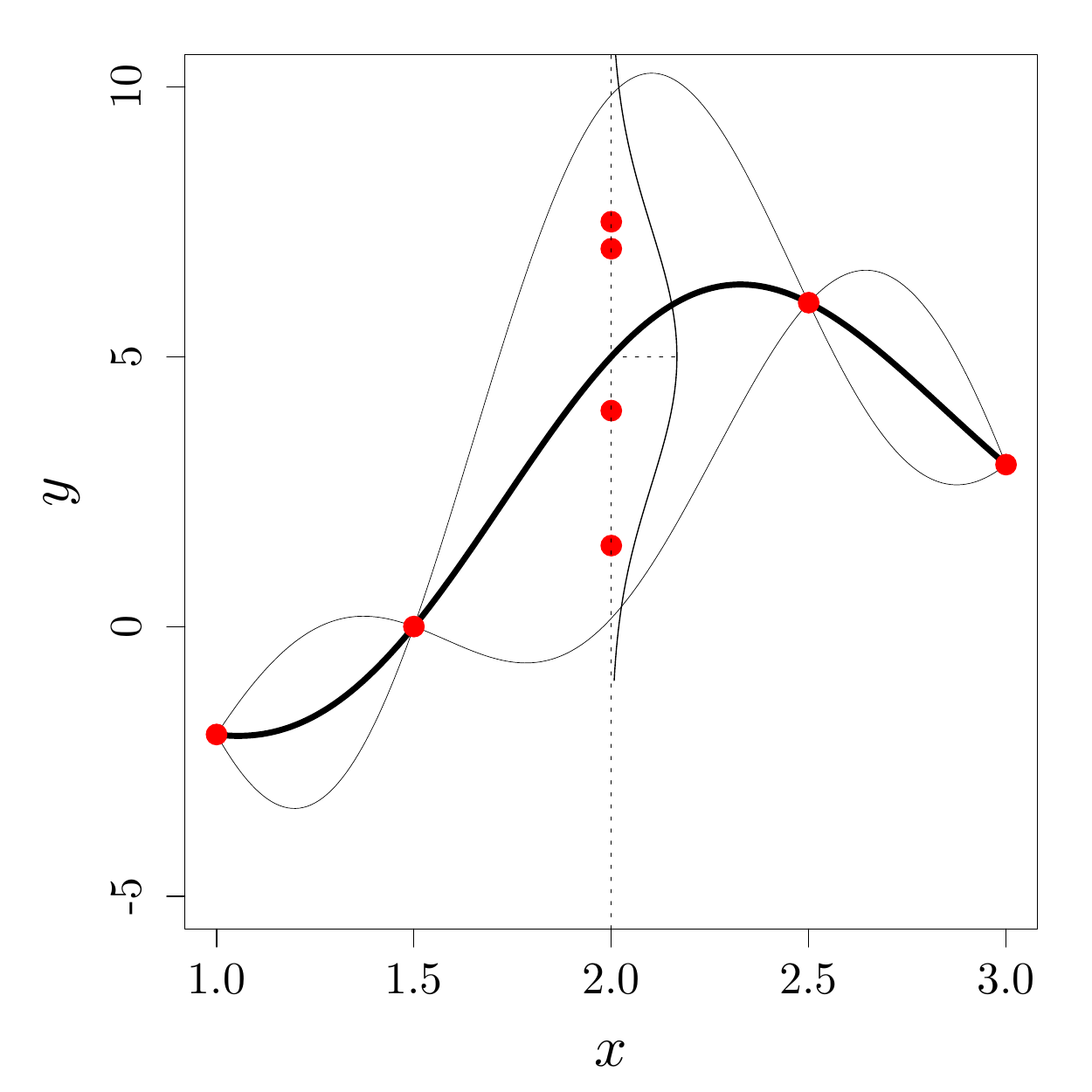}
   \caption{
Distribution-wise GP,  $m^{Dist}(x)$ (thick line) $\pm 2 \sqrt{v^{Dist}(x)}$ (thin lines). 
At the redundant point $x=2$, the outputs are 1.5, 4, 7 and 7.5.
The mean of the distribution-wise GP passes through the average of outputs.
Contrarily to PI (cf. Figure~\ref{PIaveraging}), distribution-wise GP preserves the empirical variance: 
the kriging variance at $x=2$ is equal to $\overline{s_{x=2}^2}=5.87$. \label{ideal_model}
}
\end{figure}

An example of distribution-wise GP is plotted in Figure \ref{ideal_model} where the site outputs empirical means and variances are used in the model. 

So far, we have observed that both $v^{Dist}$ and $v^{Nug}$ are non-zero at repeated points. However, there is a fundamental difference between the behaviors of a distribution-wise GP and a GP regularized by nugget; as the number of observations $N_i$ at a redundant point $\textbf{x}^i$ increases, $v^{Nug}(\textbf{x}^i)$ tends to 0 while $v^{Dist}(\textbf{x}^i)$ remains equal to $\sigma^2_{Z_i}$.   

\begin{figure}[htpb]
\centering
\includegraphics[width=0.49\textwidth]{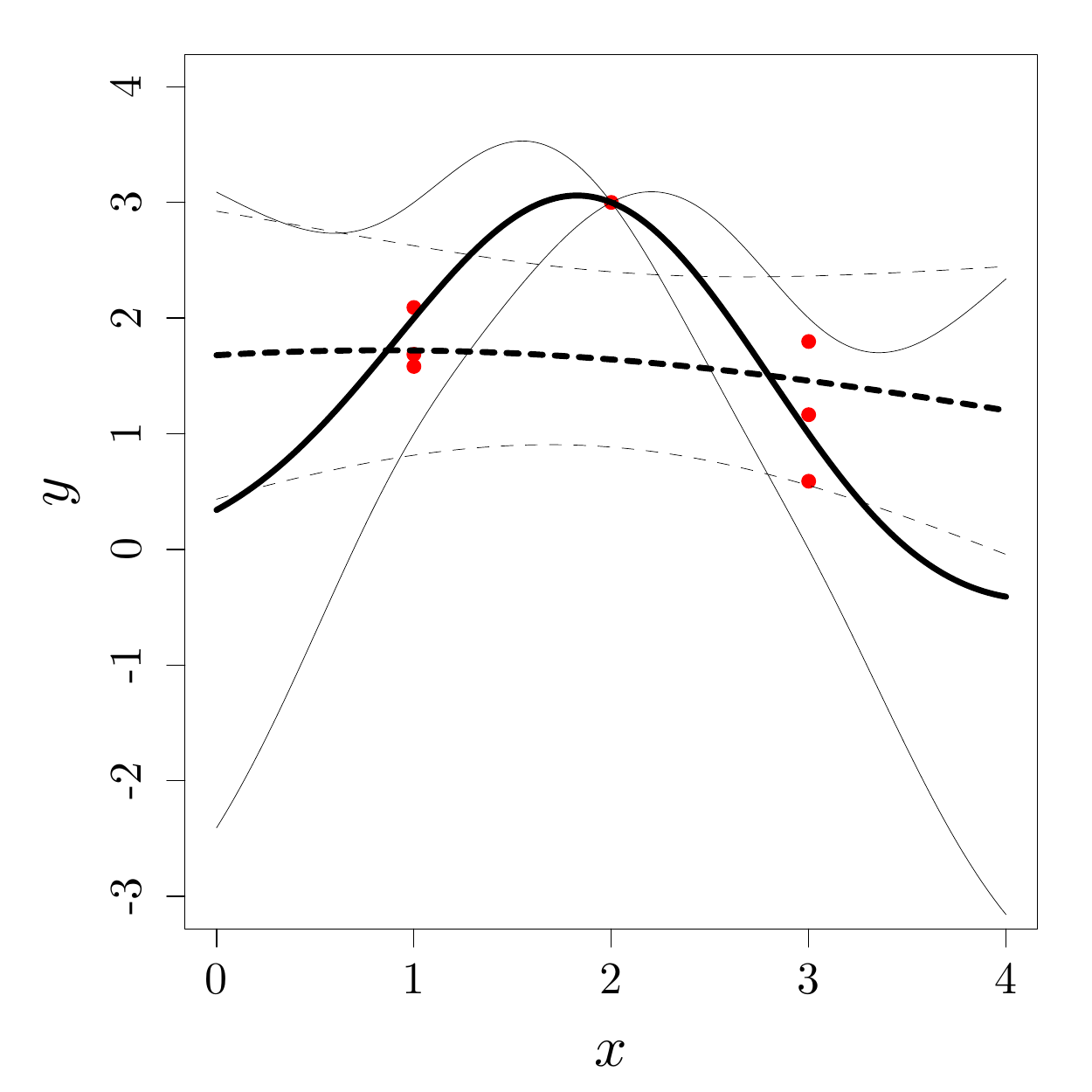}
\includegraphics[width=0.49\textwidth]{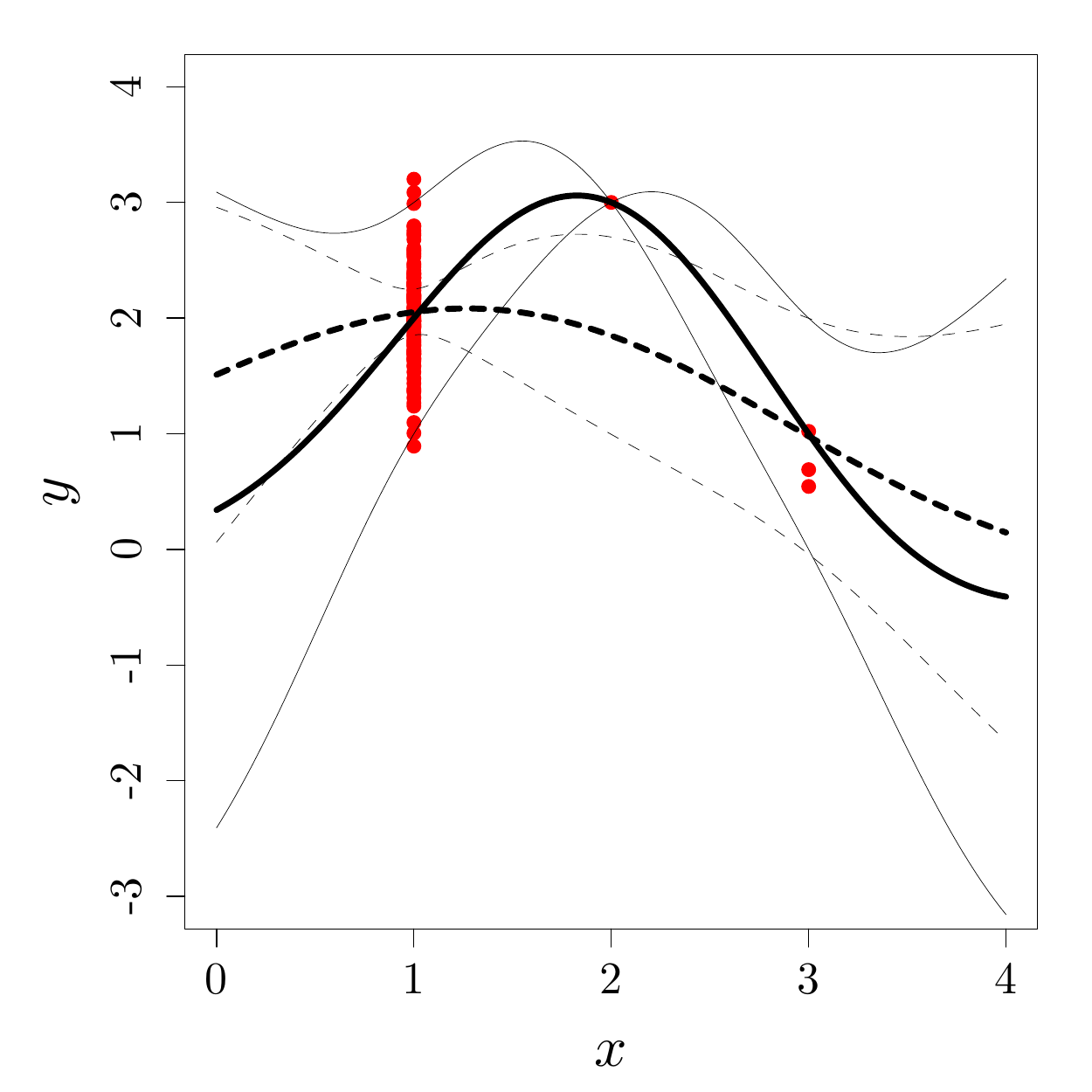}
\caption{Distribution-wise GP (solid lines) versus a GP model regularized by nugget (dashed lines). At $x=1$, the number of repeated points is 3 on the left and 100 on the right. $v^{Nug}(x=1)$ (thin dashed lines) shrinks as the number of repeated points increases while $v^{Dist}(x=1)$ remains constant.}
\label{DW_noise}
\end{figure}

This can be analytically seen by assuming that there is only one location site, $\textbf{x}^1$, with several observations, say $n$. In this situation, the correlation between every two observations is one and so, the kriging variance regularized by nugget at $\textbf{x}^1$ is 
\begin{equation}
v^{Nug}(\textbf{x}^1) =  \sigma^2 \left( 1 - [1, \ldots, 1]\left(\textbf{R} + \tau^2/\sigma^2\textbf{I}\right)^{-1} [1, \ldots, 1]^\top \right).
\label{V_nug_1}
\end{equation}
Here, the correlation matrix $\textbf{R}$ is a matrix of 1's with only one strictly positive eigenvalue equal to $\lambda_1 = n$, all other 
eigenvalues being equal to 0. 
The eigenvector associated to $\lambda_1$ is $(1, \dots, 1)^\top / \sqrt{n}$. Adding nugget will increase all the eigenvalues of $\textbf{R}$ by $\tau^2/\sigma^2$. 
In Equation (\ref{V_nug_1}) one can replace $(\textbf{R} + \tau^2/\sigma^2\textbf{I})^{-1}$ by its eigendecomposition that is,
\begin{equation}
 \begin{bmatrix}
  1/\sqrt{n} &  \\ 
  \vdots & \textbf{W} \\ 
  1/\sqrt{n} &  \\ 
  \end{bmatrix}  
  \begin{bmatrix}
  \sigma^2/n\sigma^2+\tau^2 &  & & \textbf{0}\\ 
   & \sigma^2/\tau^2 & & \\ 
   & & \ddots & \\
  \textbf{0} &  & & \sigma^2/\tau^2 \\ 
  \end{bmatrix}
  \begin{bmatrix}
  1/\sqrt{n} & \hdots & 1/\sqrt{n} \\ 
  & \textbf{W}^\top &  \\ 
  \end{bmatrix}~,
\end{equation}     
which yields 
\begin{equation}
v^{Nug}(\textbf{x}^1) =  \frac{\tau^2}{n\sigma^2 + \tau^2} \sigma^2,
\end{equation}
since $[1, \ldots, 1]$ is perpendicular to any of the other eigenvectors making the columns of $\textbf{W}$. Consequently, $v^{Nug}(\textbf{x}^1) \rightarrow 0$ when $n \rightarrow \infty$.
Figure \ref{DW_noise} further illustrates the difference between distribution-wise and nugget regularization models in GPs. 
The red bullets are data points generated by sampling from the given distribution of $\textbf{Z}$'s, 
\begin{equation*}
\textbf{Z} \quad \sim  \quad \mathcal N \left(
\begin{bmatrix} 
2 \\ 3 \\ 1 \\
\end{bmatrix} 
~ , ~
\begin{bmatrix} 
0.25 & 0 & 0\\
0 & 0 & 0 \\
0 & 0 & 0.25 \\
\end{bmatrix} 
\right)
\end{equation*}
and the right plot has more data points at $x=1$ than the left plot.
We observe that the distribution-wise GP model is independent from the number of data points and, in that sense, it ``interpolates the distributions'':
the conditional variance of the distribution-wise GP model does not change with the increasing number 
of data points at $x=1$ while the variance of the GP model regularized by nugget decreases; the mean of the distribution-wise GP is the same on the left and right plots but that of the GP regularized by nugget changes and tends to the mean of the distribution as the number of data points grows.

\section{Conclusions} \label{conclusion}
This paper provides a new algebraic comparison of pseudoinverse and nugget regularizations, the 
two main methodologies to overcome the degeneracy of the covariance matrix in Gaussian processes (GPs). 
By looking at non invertible covariance matrices as the limit for ill-conditionned covariance matrices, 
we have defined redundant points. Clear differences between pseudoinverse and nugget regularizations
have arised: contrarily to GPs with nugget, GPs with pseudoinverse average the values of
outputs and have zero variance at redundant points; 
in GPs regularized by nugget, the discrepancy between the model and the data turns into a 
departure of the GP from observed outputs throughout the domain; 
in GPs regularized by pseudoinverse, this departure only occurs at the redundant points
and the variance is zero there.
Some guidelines have been given for choosing a regularization strategy.

In the last part of the paper, we have proposed a new regularization strategy for GPs, 
the distribution-wise GP. This model interpolates normal distributions instead 
of data points. It does not have the drawbacks of nugget and pseudoinverse regularizations: 
it not only averages the outputs at repeated points but it also preserves 
the repeated points variances.

Distribution-wise GPs shed a new light on regularization, which starts with the creation 
of repeated points by clustering. 
A potential benefit is the reduction in covariance matrix size.  
Further studying distribution-wise GPs is the main continuation of this work.

\section*{Acknowledgments}
The authors would like to acknowledge support by the French national research agency (ANR) within the Modèles Numérique project ``NumBBO- Analysis, Improvement and Evaluation of Numerical Blackbox Optimizers".

\appendix


\section{Examples of redundant points}
\label{app-redundant}
This Appendix gives easily interpretable examples of DoEs with associated kernels that make the covariance matrix non-invertible.
The eigenvalues, eigenvectors and orthogonal projection matrix onto the image space (cf. also Section~\ref{sec-redundant}) are described.

\subsection{Repeated points}
\label{app-redund-repeated}
\begin{figure}[h] 
\begin{minipage}[c]{0.49\textwidth}
 \includegraphics[width=6cm]{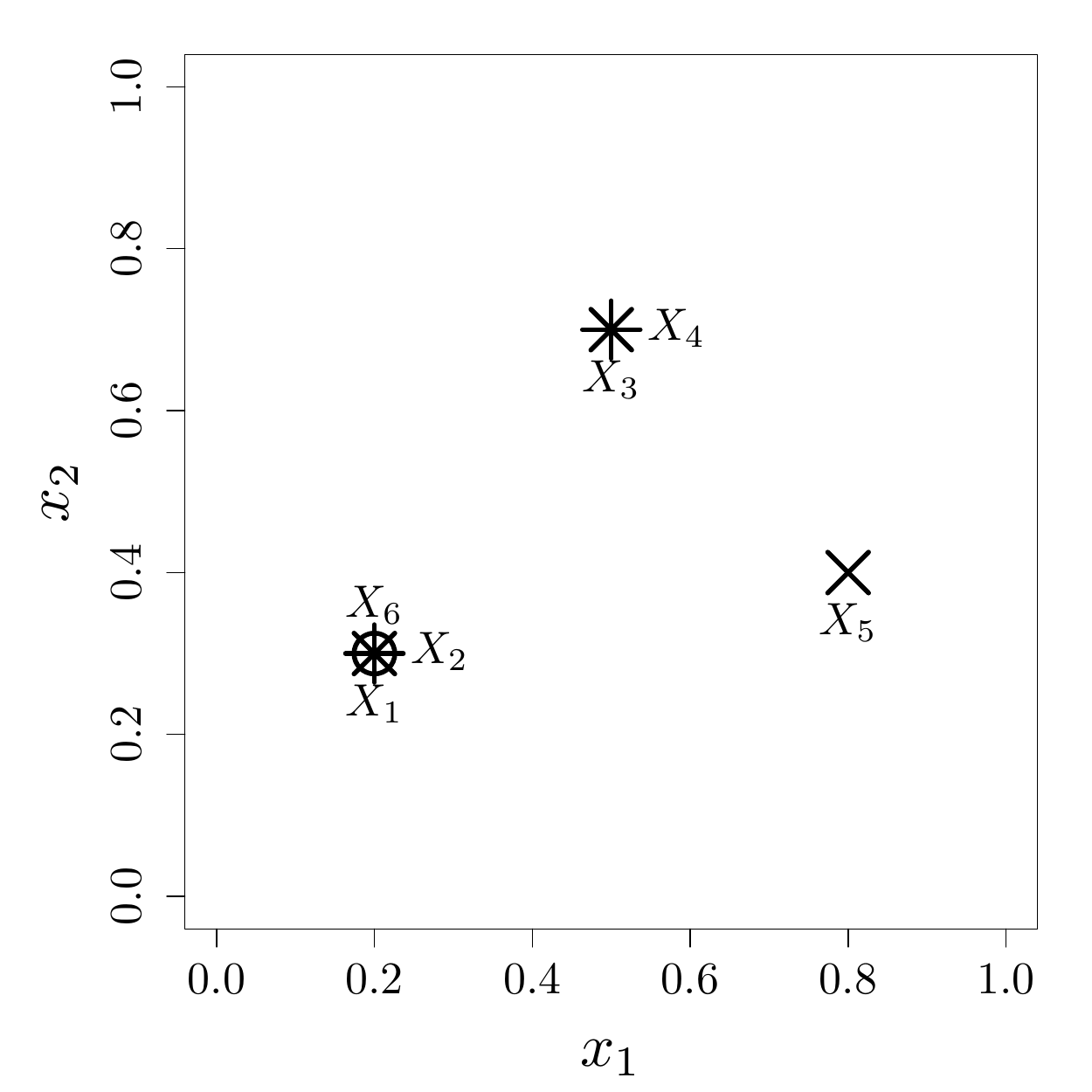} 
\end{minipage}
\begin{minipage}[c]{0.49\textwidth}
\begin{center}
$$ k(\textbf{x},\textbf{x'}) = \exp \left(-\frac{(x_1-x_1')^2}{2 \times.25^2} \right) \times \exp \left(-\frac{(x_2-x_2')^2}{2 \times.25^2} \right)$$
$
  \textbf X =
\begin{bmatrix}{}
  0.20 & 0.30 \\ 
  0.20 & 0.30 \\ 
  0.50 & 0.70 \\ 
  0.50 & 0.70 \\ 
  0.80 & 0.40 \\ 
  0.20 & 0.30 \\ 
  \end{bmatrix}
$
\end{center}
\end{minipage}
\caption{Kernel and DoE of the repeated points example
\label{fig-repeated}}
 \end{figure}

Repeated design points are the simplest example of redundancy in a DoE since columns of the covariance matrix $\textbf{c}$
are duplicated. 
An example is given in Figure~\ref{fig-repeated} with a two-dimensional design, and a classical 
squared exponential kernel.
The eigenvalues and eigenvectors of the covariance matrix associated to Figure~\ref{fig-repeated}
are
\begin{equation*}
  \bm\lambda=
\begin{bmatrix}{}
  3.12 \\ 
  1.99 \\ 
  0.90 \\ 
  0.00 \\ 
  0.00 \\ 
  0.00 \\ 
  \end{bmatrix}
~,~
  \textbf V=
\begin{bmatrix}{}
  -0.55 & 0.19 & 0.00 \\ 
  -0.55 & 0.19 & 0.00 \\ 
  -0.22 & -0.64 & -0.21 \\ 
  -0.22 & -0.64 & -0.21 \\ 
  -0.09 & -0.28 & 0.96 \\ 
  -0.55 & 0.19 & 0.00 \\ 
  \end{bmatrix}
~\text{ and }
  \textbf W=
\begin{bmatrix}{}
  0.00 & -0.30 & 0.76 \\ 
  -0.71 & 0.12 & -0.39 \\ 
  -0.04 & 0.66 & 0.26 \\ 
  0.04 & -0.66 & -0.26 \\ 
  0.00 & 0.00 & 0.00 \\ 
  0.71 & 0.18 & -0.37 \\ 
  \end{bmatrix}
~,
\end{equation*}
with the orthogonal projection matrix onto $Im(\textbf{C})$
\begin{equation*}
  \textbf{VV}^\top=
\begin{bmatrix}{}
  0.33 & 0.33 & 0.00 & 0.00 & 0.00 & 0.33 \\ 
  0.33 & 0.33 & 0.00 & 0.00 & 0.00 & 0.33 \\ 
  0.00 & 0.00 & 0.50 & 0.50 & 0.00 & 0.00 \\ 
  0.00 & 0.00 & 0.50 & 0.50 & 0.00 & 0.00 \\ 
  0.00 & 0.00 & 0.00 & 0.00 & 1.00 & 0.00 \\ 
  0.33 & 0.33 & 0.00 & 0.00 & 0.00 & 0.33 \\ 
  \end{bmatrix}
\end{equation*}
Points $\{1,2,6\}$ and $\{3,4\}$ are repeated and redundant.

\subsection{First additive example}
\label{app-redund-add1}
\begin{figure}[h]
\begin{minipage}[c]{0.49\textwidth}
 \includegraphics[width=7cm]{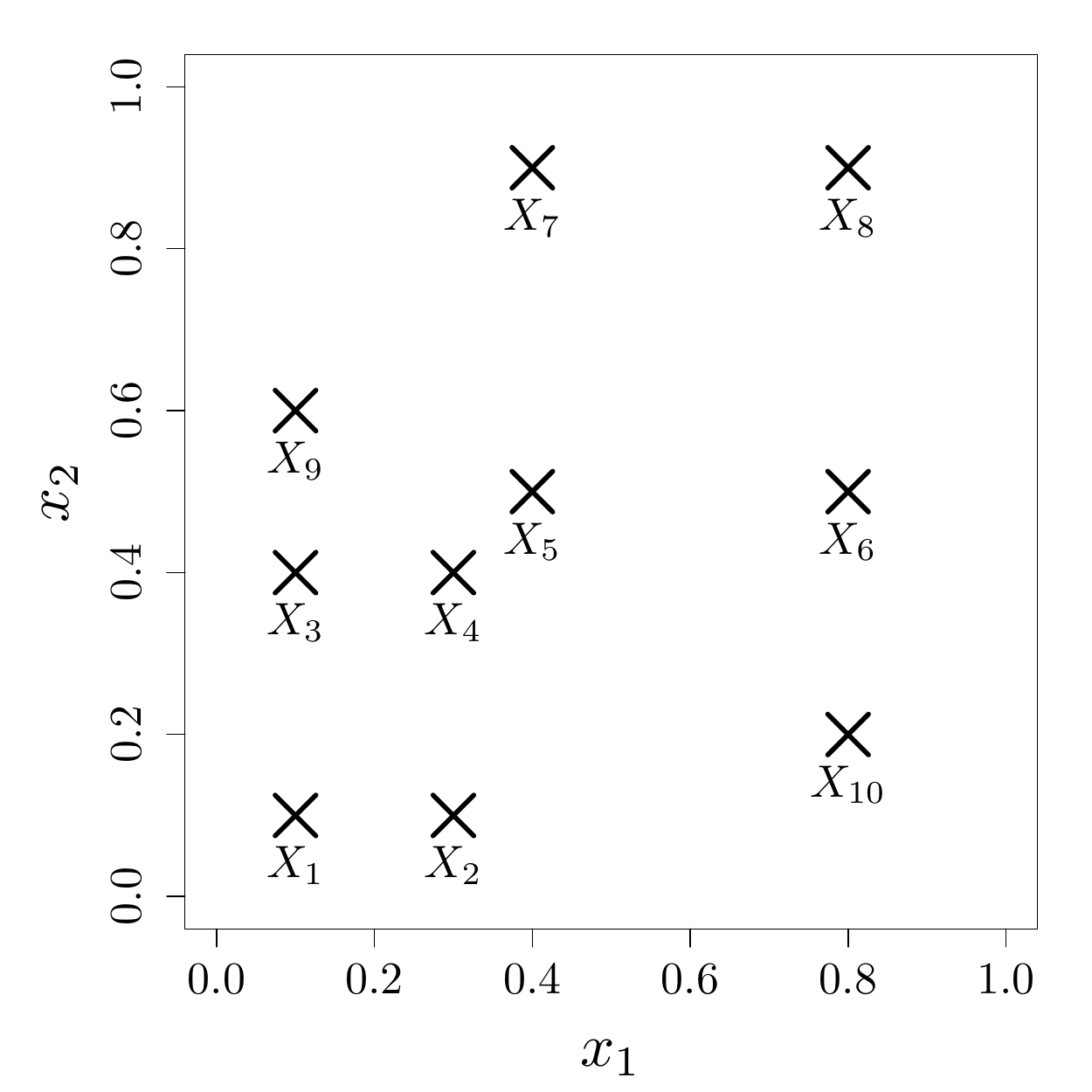}
\end{minipage}
\begin{minipage}[c]{0.49\textwidth}
\begin{center}
\begin{align}
k(\textbf{x},\textbf{x'}) = & \exp \left(-\frac{(x_1-x_1')^2}{2 \times.25^2} \right) \nonumber \\
& + \exp \left(-\frac{(x_2-x_2')^2}{2 \times.25^2} \right) \nonumber
\end{align}
$
  \textbf X =
\begin{bmatrix}{}
  0.10 & 0.10 \\ 
  0.30 & 0.10 \\ 
  0.10 & 0.40 \\ 
  0.30 & 0.40 \\ 
  0.40 & 0.50 \\ 
  0.80 & 0.50 \\ 
  0.40 & 0.90 \\ 
  0.80 & 0.90 \\ 
  0.10 & 0.60 \\ 
  0.80 & 0.20 \\ 
  \end{bmatrix}
$
\end{center}
\end{minipage}
\caption{Kernel and DoE of the first additive GP example
\label{fig-add1}}
 \end{figure}

The first example of GP with additive kernel is described in Figure~\ref{fig-add1}.
As explained in Section~\ref{sec-degeneracy}, the rectangular patterns of points $\{1,2,3,4\}$ and 
$\{5,6,7,8\}$ create linear dependencies between the columns of $\textbf{C}$.
The eigenvalues and eigenvectors of the covariance matrix are,
\begin{equation*}
  \bm \lambda=
\begin{bmatrix}{}
  9.52 \\ 
  3.58 \\ 
  2.60 \\ 
  2.31 \\ 
  1.46 \\ 
  0.39 \\ 
  0.09 \\ 
  0.06 \\ 
  0.00 \\ 
  0.00 \\ 
  \end{bmatrix}
~,~
  \textbf V =
\begin{bmatrix}{}
  -0.30 & -0.32 & 0.45 & -0.15 & 0.34 & -0.10 & 0.22 & 0.40 \\ 
  -0.33 & -0.24 & 0.29 & -0.43 & -0.22 & -0.30 & -0.43 & 0.04 \\ 
  -0.38 & -0.22 & -0.01 & 0.31 & 0.22 & 0.59 & 0.17 & 0.17 \\ 
  -0.41 & -0.14 & -0.17 & 0.04 & -0.34 & 0.40 & -0.47 & -0.19 \\ 
  -0.38 & 0.01 & -0.37 & 0.03 & -0.40 & -0.29 & 0.43 & 0.18 \\ 
  -0.28 & 0.45 & 0.03 & 0.44 & -0.13 & -0.27 & -0.15 & 0.40 \\ 
  -0.25 & 0.19 & -0.38 & -0.62 & 0.11 & 0.13 & 0.30 & -0.07 \\ 
  -0.15 & 0.64 & 0.02 & -0.22 & 0.38 & 0.15 & -0.29 & 0.15 \\ 
  -0.34 & -0.13 & -0.24 & 0.26 & 0.54 & -0.43 & -0.10 & -0.51 \\ 
  -0.25 & 0.34 & 0.59 & 0.05 & -0.22 & 0.08 & 0.35 & -0.54 \\ 
  \end{bmatrix}
\end{equation*}
\begin{equation*}
\text{ and }
  \textbf W =
\begin{bmatrix}{}
  0.00 & 0.50 \\ 
  0.00 & -0.50 \\ 
  0.00 & -0.50 \\ 
  0.00 & 0.50 \\ 
  0.50 & 0.00 \\ 
  -0.50 & 0.00 \\ 
  -0.50 & 0.00 \\ 
  0.50 & 0.00 \\ 
  0.00 & 0.00 \\ 
  0.00 & 0.00 \\ 
  \end{bmatrix}
~.
\end{equation*}
The projection matrix onto the image space is
\begin{equation*}
  \textbf{VV}^\top=
\begin{bmatrix}{}
  0.75 & 0.25 & 0.25 & -0.25 & 0.00 & 0.00 & 0.00 & 0.00 & 0.00 & 0.00 \\ 
  0.25 & 0.75 & -0.25 & 0.25 & 0.00 & 0.00 & 0.00 & 0.00 & 0.00 & 0.00 \\ 
  0.25 & -0.25 & 0.75 & 0.25 & 0.00 & 0.00 & 0.00 & 0.00 & 0.00 & 0.00 \\ 
  -0.25 & 0.25 & 0.25 & 0.75 & 0.00 & 0.00 & 0.00 & 0.00 & 0.00 & 0.00 \\ 
  0.00 & 0.00 & 0.00 & 0.00 & 0.75 & 0.25 & 0.25 & -0.25 & 0.00 & 0.00 \\ 
  0.00 & 0.00 & 0.00 & 0.00 & 0.25 & 0.75 & -0.25 & 0.25 & 0.00 & 0.00 \\ 
  0.00 & 0.00 & 0.00 & 0.00 & 0.25 & -0.25 & 0.75 & 0.25 & 0.00 & 0.00 \\ 
  0.00 & 0.00 & 0.00 & 0.00 & -0.25 & 0.25 & 0.25 & 0.75 & 0.00 & 0.00 \\ 
  0.00 & 0.00 & 0.00 & 0.00 & 0.00 & 0.00 & 0.00 & 0.00 & 1.00 & 0.00 \\ 
  0.00 & 0.00 & 0.00 & 0.00 & 0.00 & 0.00 & 0.00 & 0.00 & 0.00 & 1.00 \\ 
  \end{bmatrix}
~.
\end{equation*}
The redundancy between points 1 to 4 on the one hand, and 5 to 8 on the other hand, is readily seen on the matrix.

\subsection{Second additive example}
\label{app-redund-add2}
\begin{figure}[h]
\begin{minipage}[c]{0.49\textwidth}
 \includegraphics[width=7cm]{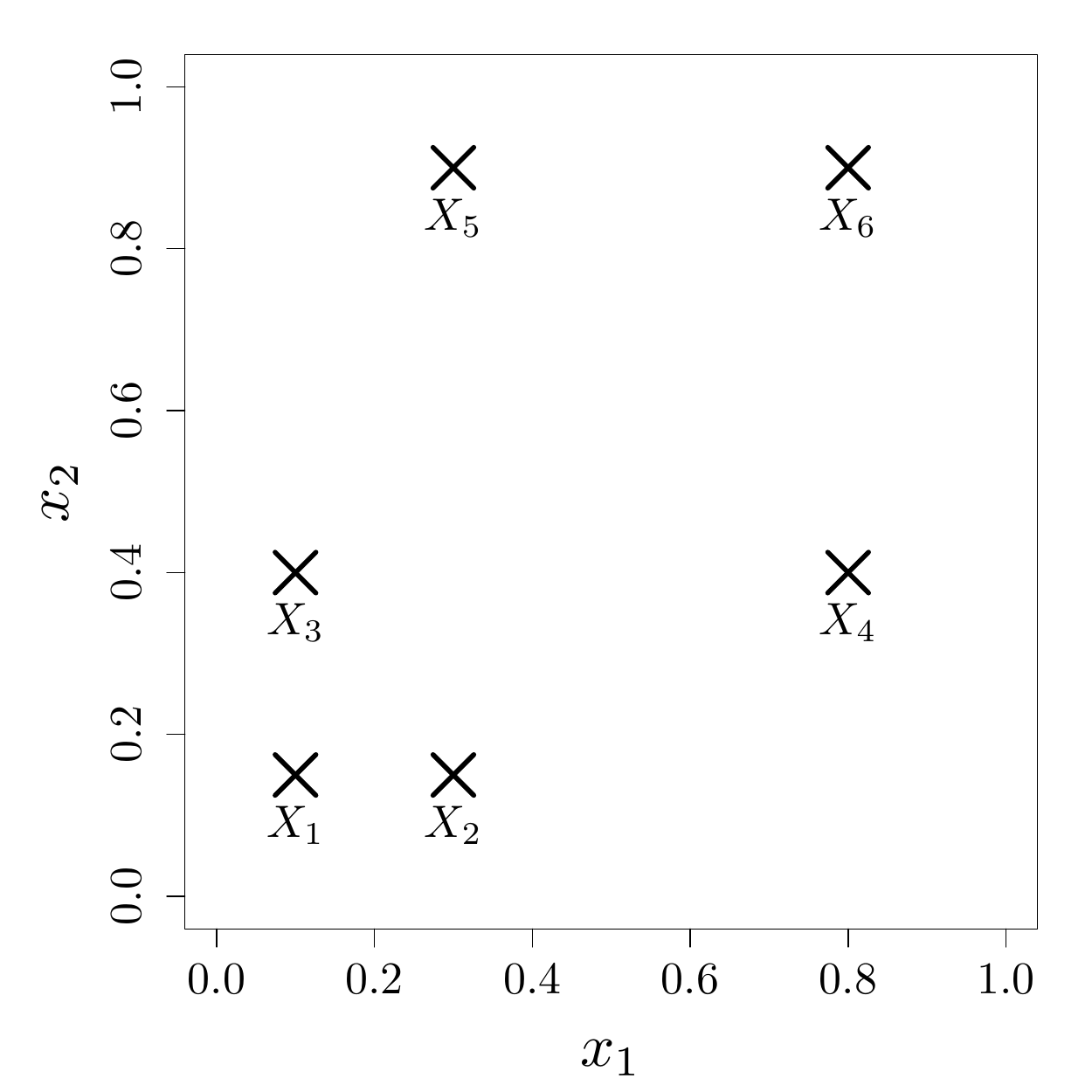}
\end{minipage}
\begin{minipage}[c]{0.49\textwidth}
\begin{center}
\begin{align}
k(\textbf{x},\textbf{x'}) = & \exp \left(-\frac{(x_1-x_1')^2}{2 \times.25^2} \right) \nonumber \\
& + \exp \left(-\frac{(x_2-x_2')^2}{2 \times.25^2} \right) \nonumber
\end{align}
$
  \textbf X =
\begin{bmatrix}{}
  0.10 & 0.15 \\ 
  0.30 & 0.15 \\ 
  0.10 & 0.40 \\ 
  0.80 & 0.40 \\ 
  0.30 & 0.90 \\ 
  0.80 & 0.90 \\ 
  \end{bmatrix}
$
\end{center}
\end{minipage}
\caption{Kernel and DoE of the second additive GP example
\label{fig-add2}}
 \end{figure}

This example shows how an incomplete rectangular pattern with additive kernels can also 
make covariance matrices singular. In Figure~\ref{fig-add2}, the point at coordinates 
$(0.3 , 0.4)$, which is not in the design, has a GP response defined twice, once by the points 
$\{1,2,3\}$ and once by the points $\{4,5,6\}$.
This redundancy in the DoE explains why $\textbf{C}$ has one null eigenvalue: 
\begin{equation*}
  \bm \lambda=
\begin{bmatrix}{}
  5.75 \\ 
  2.90 \\ 
  2.07 \\ 
  0.80 \\ 
  0.49 \\ 
  0.00 \\ 
  \end{bmatrix}
~,~
  \textbf V =
\begin{bmatrix}{}
  -0.50 & 0.34 & -0.01 & 0.18 & 0.66 \\ 
  -0.49 & 0.25 & 0.20 & 0.57 & -0.40 \\ 
  -0.48 & 0.17 & -0.29 & -0.69 & -0.01 \\ 
  -0.32 & -0.39 & -0.65 & 0.17 & -0.35 \\ 
  -0.36 & -0.28 & 0.66 & -0.33 & -0.28 \\ 
  -0.20 & -0.75 & 0.09 & 0.15 & 0.45 \\ 
  \end{bmatrix}
~,~
  \textbf W =
\begin{bmatrix}{}
  -0.41 \\ 
  0.41 \\ 
  0.41 \\ 
  -0.41 \\ 
  -0.41 \\ 
  0.41 \\ 
  \end{bmatrix}
~.
\end{equation*}
The orthogonal projection matrix onto the image space of $\textbf{C}$ tells us that 
all the points in the design are redundant,
\begin{equation*}
  \textbf{VV}^\top=
\begin{bmatrix}{}
  0.83 & 0.17 & 0.17 & -0.17 & -0.17 & 0.17 \\ 
  0.17 & 0.83 & -0.17 & 0.17 & 0.17 & -0.17 \\ 
  0.17 & -0.17 & 0.83 & 0.17 & 0.17 & -0.17 \\ 
  -0.17 & 0.17 & 0.17 & 0.83 & -0.17 & 0.17 \\ 
  -0.17 & 0.17 & 0.17 & -0.17 & 0.83 & 0.17 \\ 
  0.17 & -0.17 & -0.17 & 0.17 & 0.17 & 0.83 \\ 
  \end{bmatrix}
~.
\end{equation*}

\subsection{Periodic example}
\label{sec-periodic}
\begin{figure}[h] 
\begin{minipage}[c]{0.49\textwidth}
 \includegraphics[width=6cm]{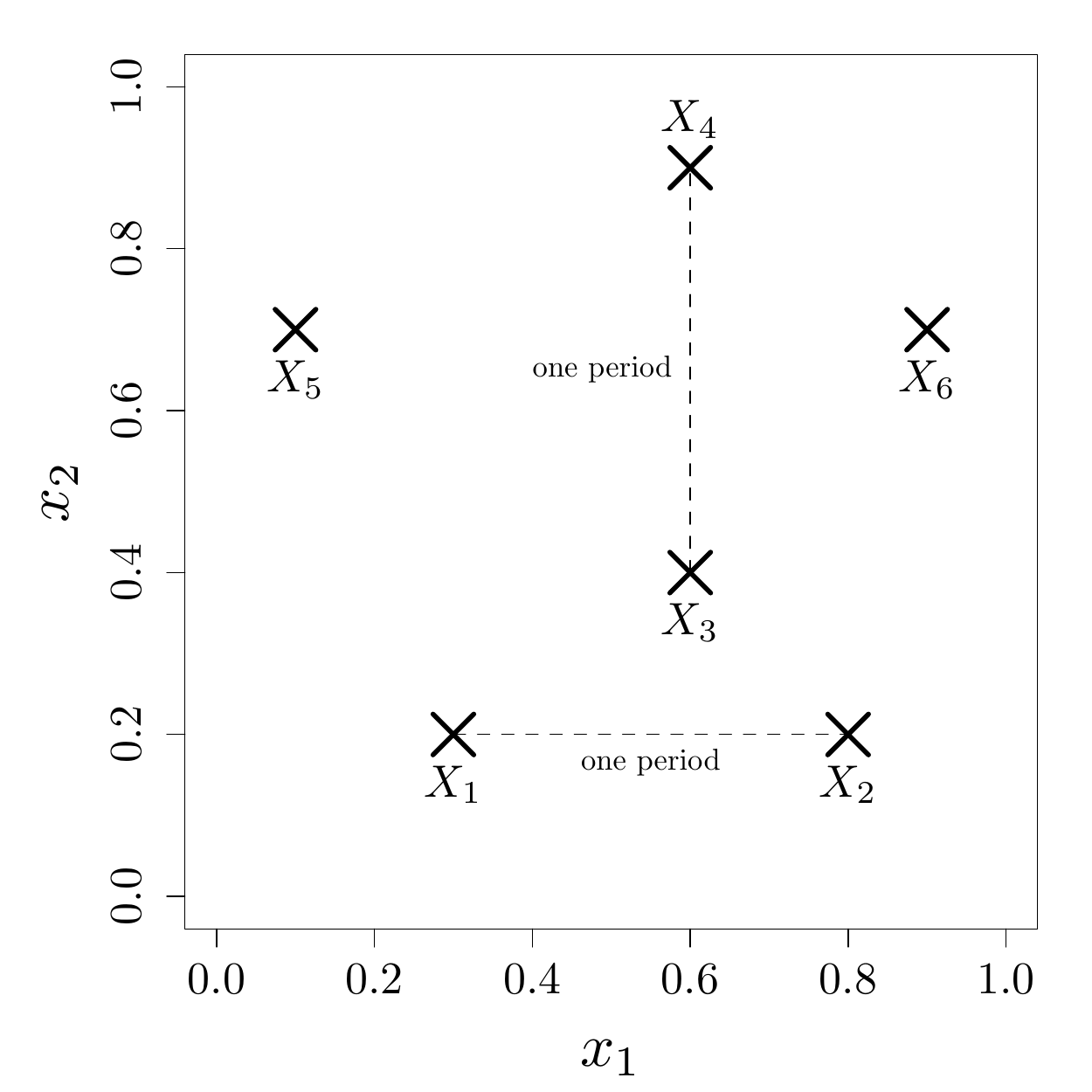} 
\end{minipage}
\begin{minipage}[c]{0.49\textwidth}
\begin{center}
\begin{eqnarray} 
k(\textbf{x},\textbf{x'}) &=& \exp \left(-\frac{\sin(4\pi (x_1-x_1'))^2}{2 \times.25^2} \right) \times \nonumber \\
 & & \exp \left(-\frac{\sin(4\pi (x_2-x_2'))^2}{2 \times.25^2} \right) \nonumber
\end{eqnarray}
$
  \textbf X =
\begin{bmatrix}{}
  0.30 & 0.20 \\ 
  0.80 & 0.20 \\ 
  0.60 & 0.40 \\ 
  0.60 & 0.90 \\ 
  0.10 & 0.70 \\ 
  0.90 & 0.70 \\ 
  \end{bmatrix}
$
\end{center}
\end{minipage}
\caption{Kernel and DoE of the periodic example
\label{fig-periodic}}
 \end{figure} 
The kernel and DoE of the periodic example are given in Figure~\ref{fig-periodic}.

The eigenvalues and eigenvectors of the associated covariance matrix $\textbf{C}$ are,
\begin{equation*}
  {\bm\lambda}= 
\begin{pmatrix}{}  
  2.00 \\ 
  2.00 \\ 
  1.01 \\ 
  0.99 \\ 
  0.00 \\ 
  0.00 \\ 
\end{pmatrix}
~ , ~
  \textbf V =
\begin{bmatrix}{}
  -0.50 & 0.50 & 0.01 & -0.01 \\ 
  -0.50 & 0.50 & 0.01 & -0.01 \\ 
  -0.50 & -0.50 & 0.01 & -0.01 \\ 
  -0.50 & -0.50 & 0.01 & -0.01 \\ 
  -0.03 & 0.00 & -0.70 & 0.72 \\ 
  0.00 & 0.00 & -0.72 & -0.70 \\ 
  \end{bmatrix}
 ~ \text{and} ~
  \textbf W =
\begin{bmatrix}{}
  0.00 & 0.71 \\ 
  0.00 & -0.71 \\ 
  0.71 & 0.00 \\ 
  -0.71 & 0.00 \\ 
  0.00 & 0.00 \\ 
  0.00 & 0.00 \\ 
  \end{bmatrix}
~.
\end{equation*} 
There are two null eigenvalues. The projector onto the image space is
\begin{equation*}
  \textbf{VV}^\top=
\begin{bmatrix}{}
  0.50 & 0.50 & 0.00 & 0.00 & 0.00 & 0.00 \\ 
  0.50 & 0.50 & 0.00 & 0.00 & 0.00 & 0.00 \\ 
  0.00 & 0.00 & 0.50 & 0.50 & 0.00 & 0.00 \\ 
  0.00 & 0.00 & 0.50 & 0.50 & 0.00 & 0.00 \\ 
  0.00 & 0.00 & 0.00 & 0.00 & 1.00 & 0.00 \\ 
  0.00 & 0.00 & 0.00 & 0.00 & 0.00 & 1.00 \\ 
  \end{bmatrix} 
\end{equation*}
which shows that points 1 and 2, on the one hand, and points 3 and 4, on the other hand, are redundant.

\subsection{Dot product kernel example}
\label{sec-linear}
The non-stationary dot product or linear kernel is \quad
$ k(\textbf{x},\textbf{x'}) = 1 + \textbf{x}^\top \textbf{x'}$. \\ \noindent
We consider a set of three one dimensional, non-overlapping, observation points:
$
  \textbf X =
\begin{bmatrix}{}
  0.20 \\ 
  0.60 \\ 
  0.80 \\ 
  \end{bmatrix}
$. \qquad The associated eigenvalues and eigenvectors are,
\begin{equation*}
  \bm \lambda=
\begin{bmatrix}{}
  3.90 \\ 
  0.14 \\ 
  0.00 \\ 
  \end{bmatrix}
~ , ~
  \textbf V=
\begin{bmatrix}{}
  -0.49 & 0.83 \\ 
  -0.59 & -0.09 \\ 
  -0.64 & -0.55 \\ 
  \end{bmatrix}
\text{~ and ~}
  \textbf W=
\begin{bmatrix}{}
  0.27 \\ 
  -0.80 \\ 
  0.53 \\ 
  \end{bmatrix}
\end{equation*}
The projection matrix onto the image space of $\textbf C$ is
\begin{equation*}
  \textbf{VV}^\top=
\begin{bmatrix}{}
  0.93 & 0.21 & -0.14 \\ 
  0.21 & 0.36 & 0.43 \\ 
  -0.14 & 0.43 & 0.71 \\ 
  \end{bmatrix}
\end{equation*}
Because there are 3 data points which is larger than $d+1=2$, all points are redundant.
With less than 3 data points, the null space of $\textbf C$ is empty.

\section{Proof of Theorem 1}
\label{sec-appnugget}

Before starting the proof, we need equations resulting from the positive definiteness 
of the covariance matrix $\textbf{C}$:
\begin{align}
&\textbf{y}=\textbf{P}_{Null(\textbf{C})}\textbf{y} + \textbf{P}_{Im(\textbf{C})}\textbf{y}\\
&\textbf{P}_{Im(\textbf{C})}\textbf{y}=\sum \limits_{i=1}^{n-N+k}\langle \textbf{y}, \textbf{V}^i \rangle \textbf{V}^i \\
&\textbf{P}_{Null(\textbf{C})}\textbf{y}=\sum \limits_{i=1}^{N-k}\langle \textbf{y}, \textbf{W}^i \rangle \textbf{W}^i \\
&\left\Vert \textbf{P}_{Null(\textbf{C})}\textbf{y}\right \Vert ^2 = \left\Vert \textbf{y} - \textbf{P}_{Im(\textbf{C})}\textbf{y} \right\Vert ^2, 
\end{align}
\noindent where $\langle.,.\rangle$ denotes the inner product.

The natural logarithm of the likelihood function is
\begin{eqnarray} 
\ln L(\textbf{y}\vert \bm{\theta}, \sigma^2)=-\frac{n}{2}\ln(2\pi) - \frac{1}{2}\ln\vert \textbf{C} \vert - \frac{1}{2}\textbf{y}^\top \textbf{C}^{-1}\textbf{y},
\end{eqnarray} 
where after removing fixed terms and incorporating nugget effect, becomes:
\begin{eqnarray} \label{E:nugget_likelihood}
-2\ln L(\textbf{y}\vert \tau^2) \approx \ln \left(\left\vert \textbf{C}+\tau^2 \textbf{I}\right\vert \right) + \textbf{y}^\top \left(\textbf{C}+\tau^2 \textbf{I}\right)^{-1}\textbf{y}.
\end{eqnarray} 

The eigenvalue decomposition of matrix $\textbf{C}+\tau^2 \textbf{I}$ in (\ref{E:nugget_likelihood}) consists of
\begin{align}
&\left(\textbf{V}^1, ..., \textbf{V}^{n-N+k}, \textbf{W}^{1}, ..., \textbf{W}^{N-k} \right)\\
&\bm{\Sigma}=diag(\tau^2 + \lambda_{1}, ..., \tau^2 + \lambda_{n-N+k}, \displaystyle\underbrace{\tau^2, ..., \tau^2}_{N-k}).
\end{align}
If Equation~(\ref{E:nugget_likelihood}) is written based on the eigenvalue decomposition, we have
\begin{equation}
-2\ln L(\textbf{y}\vert \tau^2) \approx  \sum \limits_{i=1}^n \ln(\tau^2 + \lambda_i) + \frac{1}{\tau^2} \sum \limits_{i=1}^{N-k} \langle \textbf{y}, \textbf{W}^i \rangle ^2 + \sum\limits_{i=1}^{n-N+k}\frac{\langle\textbf{y},\textbf{V}^i \rangle^2}{\tau^2+\lambda_i},
\end{equation}
\noindent or equivalently
\begin{eqnarray}\label{likelihood_with_nugget_2}
-2\ln L(\textbf{y}\vert \tau^2) \approx  \sum \limits_{i=1}^n \ln(\tau^2 + \lambda_i) + \frac{1}{\tau^2} \left\Vert \textbf{y} - \textbf{P}_{Im(\textbf{C})}\textbf{y}\right\Vert ^2 + \sum\limits_{i=1}^{n-N+k}\frac{\langle\textbf{P}_{Im(\textbf{C})}\textbf{y},\textbf{V}^i \rangle^2}{\tau^2+\lambda_i},
\end{eqnarray}
\noindent with the convention $\lambda_{n-N+k+1}=\lambda_{n-N+k+2}=...=\lambda_{n}=0$. In the above equations, $\approx$ means ``equal up to a constant''. 
Based on (\ref{projection_IM}), the term $\textbf{y} - \textbf{P}_{Im(\textbf{C})}\textbf{y}$ in Equation (\ref{likelihood_with_nugget_2}) is
 \begin{eqnarray}\label{projection_Nul}
\textbf{y} - \textbf{P}_{Im(\textbf{C})}\textbf{y}=
\begin{bmatrix}
y_1 - \overline{y}_1 \\
\vdots \\
y_{N_1} - \overline{y}_1 \\
\vdots \\
y_{N_1+...+N_{k-1}+1} - \overline{y}_k \\
\vdots \\
y_{N_1+...+N_k} - \overline{y}_k \\
0 \\
\vdots \\
0
\end{bmatrix},
\end{eqnarray}
\noindent where $\overline{y}^i,\ i=1,..., k,$ designates the mean of response values at location $i$.

\noindent According to Equations (\ref{projection_Nul}) and (\ref{var_redundant}), $\left\Vert \textbf{y} - \textbf{P}_{Im(\textbf{C})}\textbf{y}\right\Vert ^2~=~\sum\limits_{i=1}^k N_is_i^2$. Hence, Equation~(\ref{likelihood_with_nugget_2}) using $s_i^2$ is updated as
\begin{eqnarray}
-2\ln L(\textbf{y}\vert \tau^2) \approx  \sum \limits_{i=1}^n \ln(\tau^2 + \lambda_i) + \frac{1}{\tau^2} \sum\limits_{i=1}^k N_is_i^2  + \sum\limits_{i=1}^{n-N+k}\frac{\langle\textbf{P}_{Im(\textbf{C})}\textbf{y},\textbf{V}^i \rangle^2}{\tau^2+\lambda_i}.
\end{eqnarray}

Let function $\Delta(\tau^2)$ express the difference between$-2\ln L(\textbf{y}\vert\tau^2)$ and $-2\ln L(\textbf{y}^+\vert \tau^2)$. 
Remark that $\textbf{P}_{Im(\textbf{C})}\textbf{y}= \textbf{P}_{Im(\textbf{C})}\textbf{y}^+$ because of our hypothesis $\overline{y}^i={\overline{y}^+}{}^i~, i=1, ..., k$. The function $\Delta(\tau^2)$ is defined as 
\begin{equation}
\Delta(\tau^2)~\equiv~-2\ln L(\textbf{y}^+\vert \tau^2)+ 2\ln L(\textbf{y}\vert \tau^2)~=~\frac{1}{\tau^2}\sum\limits_{i=1}^{k}N_i\left({s_i^+}^2 - s_i^2\right),
\end{equation}
and is monotonically decreasing.

Now we show that $\widehat{\tau^+}^2$, the ML estimation of nugget from $\textbf{y}^+$, is never smaller than $\hat{\tau}^2$, the ML estimation of nugget from $\textbf{y}$.
Firstly, $\widehat{\tau^+}^2$ cannot be smaller than $\hat{\tau}^2$.  Indeed, if $\tau^2 \leq \hat{\tau}^2$, then
\begin{align} \label{last_proof}
-2\ln L(\textbf{y}^+\vert\tau^2)&=-2\ln L(\textbf{y}\vert\tau^2)+\Delta(\tau^2) \\
\nonumber&\geq -2\ln L(\textbf{y}\vert \hat{\tau}^2)+\Delta(\tau^2) \\
\nonumber&\geq -2\ln L(\textbf{y}\vert \hat{\tau^2})+\Delta(\hat{\tau}^2) \\
\nonumber&=-2\ln L(\textbf{y}^+\vert\hat{\tau}^2), 
\end{align}
which shows that $\widehat{\tau^+}^2 \geq \hat{\tau}^2$.
Secondly, if ${s_i^+}^2$ is strictly larger than $s_i^2$, then  $\widehat{\tau^+}^2 > \hat{\tau}^2$ because the slope of $-2\ln L(\textbf{y}^+\vert \tau^2)$ is strictly negative at $\tau^2=\hat{\tau}^2$: The derivative of $-2\ln L(\textbf{y}^+\vert \tau^2)$ with respect to $\tau^2$ can be written as
\begin{align}
\frac{d}{d\tau^2}\left( -2\ln L(\textbf{y}^+\vert \tau^2) \right)=\frac{d}{d\tau^2}\left(-2\ln L(\textbf{y} \vert \tau^2)\right) + \frac{d\Delta(\tau^2)}{d\tau^2}.
\end{align}
Since $\hat{\tau}^2=\arg\min -2\ln L(\textbf{y}\vert \tau^2)$, the second term in the right hand side of the above equation is equal to zero. Therefore, the derivative of $-2\ln L(\textbf{y}^+\vert \tau^2)$ with respect to $\tau^2$ reduces to 
\begin{eqnarray}
\frac{d}{d\tau^2}\left( -2\ln L(\textbf{y}^+\vert \hat{\tau}^2) \right)=\frac{d}{d\tau^2}\left(\frac{1}{\tau^2}\sum\limits_{i=1}^{k}N_i\left({s_i^+}^2 - s_i^2\right)\right)~=~\frac{-1}{\tau^4}\sum\limits_{i=1}^{k}N_i\left({s_i^+}^2 - s_i^2\right).
\end{eqnarray}
The above derivative is strictly negative because ${s_i^+}^2 - s_i^2$ is positive and the proof is complete. $\square$

\bibliography{biblio}
\bibliographystyle{plain}
\end{document}